\documentclass{article}

\usepackage{amsmath, amsthm, amssymb}
\usepackage{xfrac}

\usepackage{marginnote}
\usepackage{graphicx}

\usepackage[table]{xcolor} 

\usepackage[T1]{fontenc}

\usepackage{fullpage}

\usepackage[numbers]{natbib}
\usepackage{algorithmic}
\usepackage{algorithm}
\usepackage{url}

\usepackage{soul}

\usepackage{hyperref}

\usepackage{xspace}

\newtheorem{definition}{Definition}
\newtheorem{prop}{Proposition}
\newtheorem{myth}{Theorem}
\newtheorem{myco}{Corollary}

\newcommand{\Probab}[1]{\mathcal{P}({#1})}
\newcommand{\Pcond}[2]{\Probab{{#1}\mid{#2}}}
\newcommand{\Pconj}[2]{\Probab{{#1},{#2}}}
\renewcommand{\~}[1]{\overline{#1}}

\newcommand{\hC}[1]{\cellcolor{gray!85}{#1}}
\newcommand{\mC}[1]{\cellcolor{gray!40}{#1}}
\newcommand{\lC}[1]{\cellcolor{gray!20}{#1}}

\newcommand{\PR}{\textsc{pr}}

\newcommand{\cause}{causation\xspace}
\newcommand{\caprese}

\newcommand{\New}[1]{{\color{black}{#1}}} 
\newcommand{\NewLast}[1]{{\color{black}{#1}}}

\begin{document}

\title{Inferring tree causal models of cancer progression\\ with
  \textcolor{black}{probability raising} \\
}


\author{Loes Olde Loohuis\footnote{Equal contributors.} \footnote{L. Olde Loohuis is currently at the Center for Neurobehavioral Genetics,  University of California, Los Angeles, USA.} \\
Department of Computer Science \\
City University New York, The Graduate Center\\
New York, USA.
\and 
Giulio Caravagna$^*$ $\quad$ Alex Graudenzi$^*$\\
Daniele Ramazzotti$^*$ $\quad$
Giancarlo Mauri\\
Marco Antoniotti\\
Dipartimento di Informatica Sistemistica e Comunicazione\\
Universit\`a degli Studi Milano-Bicocca\\
Milano, Italy.
\and
Bud Mishra \\
Courant Institute of Mathematical Sciences \\
New York University\\
New York, USA.
}

\date{}

\maketitle


\begin{abstract}
Existing techniques to reconstruct tree models of progression for
accumulative processes, such as cancer, seek to estimate \cause 
by combining correlation and a frequentist notion of temporal
priority. In this paper, we define a novel theoretical framework
\New{called CAPRESE (CAncer PRogression Extraction with Single Edges)} 
to reconstruct such models based on the notion of probabilistic \cause
defined by Suppes.  
We consider a general reconstruction setting complicated by the
presence of noise in the data due to 
biological variation, as well as experimental or measurement
errors. 
\New{To improve
  tolerance to} noise we define and use a
shrinkage-like estimator.
\New{We prove the correctness of our algorithm by showing asymptotic
  convergence to the correct tree under mild constraints on the level
  of noise. Moreover,} on synthetic data, we show that our approach
outperforms the state-of-the-art, that it is efficient even with a
relatively small number of samples and that its performance quickly converges to its asymptote as the number of samples increases. For real cancer datasets obtained with different technologies, we highlight biologically significant differences in the progressions inferred with respect to other competing techniques and we also show how to validate conjectured biological relations with progression models. 
\end{abstract}



\section*{Introduction}

Cancer is a disease of evolution. 
Its initiation and progression 
are caused by dynamic somatic alterations to the genome manifested as point 
mutations, structural alterations, DNA methylation and histone 
modification changes \cite{Hanahan_Weinberg_2011}. 

These genomic alterations are generated by random processes, and since individual tumor cells compete for space and resources,
the fittest variants are naturally selected for.
For example, if through some mutations a 
cell acquires 
the ability to ignore anti-growth signals from the
body, this cell may thrive and
divide, and its progeny may eventually dominate some part(s) of the
tumor. 
This \emph{clonal expansion} can be seen as a
\emph{discrete state} of the cancer's progression, marked by the
acquisition of a set of genetic events. 
Cancer progression can then be thought  of as a sequence
of these
discrete  steps, where the tumor acquires certain distinct properties at
each state. 
Different progression sequences are possible, but some are  more
common than others, and not every order is
viable \cite{luo_principles_2009}.


In the last two decades, many specific genes and genetic mechanisms that are involved in different types of cancer have been
identified (see e.g. \cite{vogelstein2004cancer,frank2007} for an overview of
common cancer genes and \cite{bell2011integrated, imielinski2012mapping} for specific genetic analyses of
ovarian carcinoma and lung adenocarcinoma, respectively), 
and \emph{therapies} targeting the activity of these genes are now being developed at a fast pace \cite{luo_principles_2009}. However, unfortunately, the {\em causal and temporal relations} 
among the genetic events driving cancer progression remain largely elusive. 

The main reason for this state of affairs is that information revealed in the data is usually obtained only at one (or a few)
points in time, rather than over the course of the disease. 
Extracting this dynamic information from 
the available  {\em cross-sectional} data is  challenging,
and a combination of mathematical, statistical and computational
techniques  is needed. 
\New{In recent years, several methods to extract
progression models from cross-sectional data
have been developed, starting from the seminal work on
single-path-models by Fearon and Vogelstein
\cite{vogelstein1988genetic}. In particular, different models of
oncogenetic trees were developed over the years. At the core of some
of these methods, e.g. \cite{desper1999inferring,desper2000distance},
is the use of \emph{correlation} to identify relations among genetic
events. These techniques reconstruct \emph{tree} models of progression
as independent  acyclic paths with  branches and no
confluences. Distinct models of oncogenetic trees are instead based on
{\em maximum likelihood estimation}, e.g.,
\cite{vonHeydebreck2004,beerenwinkel2005learning,Beerenwinkel2005}. More
general {\em Markov chain} models, e.g., \cite{Hjelm_2006}, describe
more flexible probabilistic networks, despite the computationally
expensive parameter estimation.   Other recent models are Conjunctive
Bayesian Networks, CBNs
\cite{beerenwinkel2007conjunctive,gerstung2009quantifying}, that
extract \emph{directed acyclic graphs}, yet imposing specific constraints on the joint occurrence  of events. Finally, in a slightly different context, temporal models were reconstructed from 
time-course gene expression data \cite{Bar-Joseph,Ramakrishnan}. 
}

In this paper we present a novel theoretical framework \New{called CAPRESE (CAncer PRogression Extraction with Single Edges)}  to
reconstruct cumulative
progressive phenomena, such as cancer progression. 
We assume the original problem setting of \cite{desper1999inferring},
and propose a new a technique to infer {\em probabilistic progression trees} from
cross-sectional data. \New{Unlike maximum likelihood estimation-based techniques, our goal is the extraction of the {\em minimal} progression model explaining the order in which mutations occur and accumulate.} The method is technology agnostic, i.e., it can be applied to dataset 
derived from all types of (epi-)genetic data such as deep exome sequencing, bisulfite 
sequencing, SNP arrays, etc., (see Results), and takes as input a set of 
pre-selected genetic events of which the presence or the absence of each event is recorded 
for each sample. 

\New{CAPRESE is} based on two main ingredients:
$1)$ rather than using \emph{correlation} to infer progression
structures, we base our technique on a notion of \emph{probabilistic 
\cause{}}, and $2)$, to increase robustness
against noise, we adopt a {\em shrinkage-like estimator} to 
measure \cause among any pair of events.  

More specifically, with respect to our first ingredient, we adopt the  notion of (prima facie) \cause proposed by 
Suppes in \cite{suppes1970probabilistic}. Its basic intuition is simple: event 
$a$ causes event  $b$ if $(i)$ $a$ occurs \emph{before} $b$ and 
$(ii)$ the occurrence of  $a$ {\em raises the probability} of 
observing $b$. 
This is a very basic notion of probabilistic \cause that in itself does not
address many of the problems associated with it (such as asymmetry,
common causes, and screening off 
\cite{sep-causation-probabilistic}), and includes \emph{spurious} as well as
\emph{genuine} causes. 
However, as it turns out, this basic notion combined with a filter
for independent progressions starting from the same root, is an excellent tool to guide progression
extraction from cross-sectional data -- one that outperforms the
commonly used correlation-based methods. 

Probabilistic \cause was  used in biomedical applications before (e.g., to find driver genes from CNV 
data in \cite{ionita_mapping_2006}, and to extract causes from biological time series data in \cite{kleinberg-thesis}), but, 
to the best of our knowledge, never to infer \emph{progression models} in the {\em absence} 
of direct temporal information. 





The extraction problem is complicated by the
presence of both false positive and false negative observations 
(see \cite{szabo2002} for a discussion on this issue based on the 
reconstruction by \cite{desper1999inferring}), 
such as the one provided by the intrinsic
variability of {biological processes} (e.g., {\em genetic 
heterogeneity}) and {\em experimental  errors}. 
This poses a problem, because while probability raising is a very
precise tool, it, by itself, is not robust enough against noise. 
Conditional on the amount of noise,  we will rely both on probabilistic \cause and on a more robust 
(but less precise) correlation-based metric in an optimal way. 
To do this we introduce our second ingredient, a {\em shrinkage-like estimator} to 
measure \cause among any pair of events. The intuition behind this
estimator, which is closely related to a shrinkage estimator from
\cite{EmpiricalBayes}, is to find the optimal balance between
probability raising on the one hand and correlation on the other,
depending on the amount of noise. 

\New{
We prove correctness of our algorithm by showing that with increasing
sample sizes, the reconstructed tree asymptotically converges to the
correct one (Theorem 3). Under mild constraints on the noise
  rates, this result holds for 
 the reconstruction problem in the presence of uniform noise as well.
}

\New{
We also study the performance of CAPRESE in more realistic settings with limited sample
sizes. Using synthetic data,} we show that \New{under these
conditions, }our algorithm outperforms the 
state-of-the-art tree reconstruction algorithm of
\cite{desper1999inferring} (see Results). 
In particular, our shrinkage-like estimator provides, on average, an
increased robustness to noise which ensures it to outperform oncotrees
\cite{desper1999inferring}.
Performance is defined in terms of {\em structural similarity}
between the reconstructed tree and the actual tree, rather than on their induced 
distribution as is done, e.g., in
\cite{beerenwinkel2005learning}. This metric is 
especially appropriate for the goal of reconstructing 
a progression model where data-likelihood fit is secondary to ``calling'' the possibly minimal set of 
causal relations.


Also, we show that 
\New{CAPRESE} works well already with a relatively low number of samples and that its performance quickly converges to its asymptote as the number of samples increases. This outcome hints at the applicability of the algorithm with relatively small 
datasets without compromising its efficiency. 

\New{We remark that further analyses on synthetic data suggests that
CAPRESE outperforms a well known bayesian probabilistic graphical
  model as well (i.e., {\em Conjunctive Bayesian Networks} \cite{beerenwinkel2007conjunctive,gerstung2009quantifying}), which was originally conceived for the reconstruction of more complex topologies, e.g. DAGs, but was proven effective in reconstructing tree topologies as well \cite{Hainke:2012hwa} (see Results).}

Finally, we apply our technique to alterations assessed with both
Comparative Genomic Hybridization and Next Generation Sequencing
techniques (see Results). In the former case, we show that the
algorithm of \cite{desper1999inferring} and 
\New{CAPRESE} highlight biologically important differences in ovarian, gastrointestinal and oral cancer, but our inferences are statistically more significant. In the latter, we validate
a recently discovered relation among two key genes involved in leukemia.


\section*{Methods}

\subsection*{Problem setting} \label{sec:problem}

The set-up of the reconstruction problem is as follows. 
Assuming that we have a set $G$ of $n$ mutations ({\em events}, in probabilistic terminology)
and $s$ samples, we represent a cross-sectional dataset as an $s \times n$ binary 
matrix in which an entry $(k,l) = 1$ if the mutation $l$ was
observed in sample $k$, and $0$ otherwise. 
The problem we  solve in this paper is to extract a set of
edges $E$ yielding  a progression \emph{tree} ${\cal
  T}=(G\cup\{\diamond\}, E, \diamond)$  from this matrix which, we remark, only implicitly provides information of progression timing. 
The root of $\cal T$ is modeled using a (special)  event  $\diamond
\not \in G$ such that {\em heterogenous progression paths} or
\emph{forests} can be reconstructed. 
More precisely, we aim at reconstructing a  \emph{rooted tree}  that satisfies: $(i)$ each 
node has at most one incoming edge, $(ii)$ the root has no incoming edges $(iii)$ 
there are no {\em cycles}. 


Each progression tree subsumes a distribution of observing a subset of
the mutations in a cancer sample that can be formalized as follows:
\begin{definition}[Tree-induced distribution] \label{def:treedistrib} Let $\cal T$ be a tree and $\alpha:E\to[0,1]$ a labeling function 
denoting the independent probability of  each edge, $\cal T$ generates a distribution where the probability 
of observing a sample with the set of alterations $G^\ast \subseteq G$ is
\begin{equation}
\Probab{G^\ast} = \prod_{e \in E'} \alpha(e) \cdot 
\mathop{
\prod_{(u,v) \in E}}_{u\in G^\ast, v\not \in G}
 \Big[1-\alpha(u,v)\Big]
\end{equation}
where all events in $G^\ast$ are assumed to be reachable from the
root $\diamond$, and $E' \subseteq E$ is the set of edges connecting the root
to the events in $G^\ast$. 
\end{definition}

We would like to emphasize two properties related to the tree-induced
distribution. 
First, the distribution subsumes that, given any oriented edge 
$(a \to b)$, an observed sample contains alteration $b$ with probability 
$\Probab{a}\Probab{b}$, that is the probability of observing $b$ after $a$. For this reason, 
if $a$ causes $b$, the probability of observing $a$ will be greater than the probability 
of observing $b$ accordingly to the temporal priority
principle which  states that all causes must precede, in time, their effects \cite{DirectionTime}. 

Second, the input dataset is a set of samples generated, ideally, from an unknown 
distribution induced by an unknown tree or forest that we aim at reconstructing. However, in some cases, 
 it could  be that no tree exists whose induced distribution generates {\em exactly} those input data. 
When this happens, the set of observed samples slightly diverges from any tree-induced distribution.
To model these situations a notion of {\em noise} can be introduced, which depends on the context in which 
data are gathered. Adding noise to the model complicates the
reconstruction problem (see Results).

\paragraph{The {\em oncotree} approach.} \label{sec:desper} 
In \cite{desper1999inferring} Desper {\em et al.} developed a method to extract progression 
trees, named \emph{``oncotrees''}, from static CNV data. \New{In
  \cite{szabo2002} Szabo {\em et al.} extended the setting of  Desper's reconstruction
  problem to account for both {\em false positives} and {\em
    negatives} in the input data. 
}In these oncotrees, nodes
represent {CNV events} and edges correspond to possible
progressions from one  event to the next. 

The reconstruction problem is exactly as described above, and 
each tree is rooted in the special event $\diamond$. The choice of which edge to include in a tree is based on the estimator 
\begin{equation}\label{eq:btree-w}
w_{a \to b}= \log{\left[ \dfrac{\Probab{a}}{\Probab{a}+\Probab{b}} \cdot \dfrac{\Pconj{a}{b}}{\Probab{a}\Probab{b}} \right]} \, ,
\end{equation}
which  assigns to each edge $a \to b$ a weight  accounting for both the relative and  joint frequencies of the events --  thus measuring 
\emph{correlation}.  The estimator is evaluated after including  $\diamond$  to each sample of the dataset. 
In this definition the rightmost term is the (symmetric) {\em likelihood ratio} for $a$ and $b$ occurring together, while the  leftmost is the  asymmetric 
\emph{temporal priority} measured by rate of occurrence. This implicit form of timing  assumes  that, if $a$ occurs \emph{more often} than $b$, then 
it likely occurs \emph{earlier}, thus satisfying  
\[
\dfrac{\Probab{a}}{\Probab{a}+\Probab{b}} > \dfrac{\Probab{b}}{\Probab{a}+\Probab{b}}\, .
\]
An oncotree is  the rooted tree whose total weight (i.e., sum of all
the weights of the edges) is maximized, and can be reconstructed in
$O(|G|^2)$ steps using  Edmond's  algorithm  \cite{edmonds1967optimum}. 
 By construction, 
the resulting graph is a proper tree rooted in $\diamond$:  each event occurs only once, {\em confluences} are absent, i.e., any event is caused by 
at most one other event. This method has been used to derive progressions for
various cancer datasets e.g.,
\cite{kainu2000somatic,huang2002genetic,radmacher2001graph}), and even
though several methods \New{that extend this framework} exists (e.g. \cite{desper2000distance,beerenwinkel2005learning,
    gerstung2009quantifying}),  to the best of our knowledge, it is currently the only method that aims to solve exactly the same problem as the one investigated in this paper and thus provide a benchmark to
compare against.

\subsection*{A probabilistic approach to \cause} \label{sec:pr} 

We briefly review the approach to 
{probabilistic \cause}, on which our method is based.
 For an extensive discussion on this topic  we refer to \cite{sep-causation-probabilistic}. 

In his seminal work \cite{suppes1970probabilistic}, Suppes proposed the following notion.

\begin{definition}[Probabilistic \cause, \cite{suppes1970probabilistic}] \label{def:praising} 
For any two events $c$ and $e$, occurring respectively at times $t_c$ and $t_e$, under the 
mild assumptions that $0 < \Probab{c}, \Probab{e} < 1$, the event $c$ is a \emph{prima facie cause} of the event $e$ if 
it occurs \emph{before} the effect and the cause \emph{raises the probability} of the effect, i.e., 
\begin{equation}
t_c < t_e \quad \text{and} \quad \Pcond{e}{c} > \Pcond{e}{\~ c} \,.
\end{equation}
\end{definition}


\New{As discussed in \cite{sep-causation-probabilistic} the above conditions are not, in general, 
sufficient to claim that event $c$ is a cause of event $e$. In fact a
prima facie cause is either \emph{genuine} or  \emph{spurious}. In the
latter case, the fact that the conditions hold in the observations is
due either to coincidence or to the presence of a certain third
\emph{confounding factor}, related both to $c$ and to $e$
\cite{suppes1970probabilistic}. Genuine causes, instead, satisfy
Definition \ref{def:praising} and are not  screened-off by any
confounding factor. However, they need not be direct causes. See
Figure  \ref{fig:prima-facie}.}


Note that we consider cross-sectional data where no information about  
$t_c$ and $t_e$  is available, so 
\New{in our reconstruction setting} we are restricted to consider
solely the {\em probability raising} (\PR{}) property,
i.e., $\Pcond{e}{c} > \Pcond{e}{\~ c}$, \New{which makes it harder to discriminate among genuine and spurious causes}. Now we review some of its properties.
\begin{prop}[Dependency]\label{prop:PRdepend} Whenever the \PR{} holds between two events $a$ and $b$, 
then the events are {\em statistically dependent} in a positive sense, i.e., 
\begin{equation}
\Pcond{b}{a} > \Pcond{b}{\~ a} \iff \Pconj{a}{b}>\Probab{a}\Probab{b}\, .
\end{equation}
\end{prop}

This and the next proposition are  well-known facts of the \PR{}; their derivation as well as the proofs of all the results we present is in the Supplementary Materials. 
Notice that the opposite implication holds as well: when 
the events $a$ and $b$ are still dependent but in a negative sense, i.e., $\Pconj{a}{b}<\Probab{a}\Probab{b}$, the 
\PR{} does not hold, i.e., $\Pcond{b}{a} < \Pcond{b}{\~a}$. 

We would like to use the asymmetry of the \PR{}  to determine whether a pair of events $a$ and $b$ satisfy 
a \cause relation so to place $a$ before $b$ in the progression tree but, unfortunately, 
the \PR{}  satisfies the following  property.
\begin{prop}[Mutual \PR{}] \label{prop:mutualPR}
$\Pcond{b}{a} > \Pcond{b}{\~ a} \iff \Pcond{a}{b} > \Pcond{a}{\~ b}
$.
\end{prop}

That is, if $a$ raises the probability of observing $b$, then $b$ raises the
probability of observing $a$ too. 

Nevertheless, in order to determine causes and effects among the genetic events, we
can use our {\em degree of confidence} in our estimate of probability raising to decide
the direction of the \cause relationship between pairs of events. 
In other words, if $a$ raises the probability of $b$ \emph{more} than the
other way around, then $a$ is a more likely cause of $b$ than $b$ of
$a$. Notice that this is sound as long as each event has {\em at most} one cause; otherwise, {\em frequent late events} 
with more than one cause, which are rather common in biological progressive phenomena,  should be treated differently. As mentioned, the \PR{} is not symmetric, and the
\emph{direction} of probability raising depends on the relative
frequencies of the events.
We make this asymmetry precise in the following
proposition.

\begin{prop}[Probability raising and temporal priority]\label{prop:PRtemporalPrior} 
For any two events $a$ and $b$ such that the probability raising  $\Pcond{a}{b} > \Pcond{a}{\~ b}$ holds, we have 
\begin{equation} 
\Probab{a}
>\Probab{b} \iff \frac{\Pcond{b}{a}}{\Pcond{b}{\~ a}} >
\frac{\Pcond{a}{b}}{\Pcond{a}{\~ b}}\, .
\end{equation} 
\end{prop}
That is,  given that the \PR{} holds between two events, 
$a$ raises the probability of $b$ \emph{more} than $b$ raises the
probability of $a$, if and only if $a$ is observed more frequently than $b$. Notice that we use  the ratio to assess the  \PR{} inequality.
The proof of this proposition is technical and
can be found in the Supplementary Materials.
From this result it follows that if we measure the timing of an event by the
rate of its occurrence (that is, $\Probab{a} > \Probab{b}$ implies that $a$ happens
{before} $b$),  this notion of \PR{} subsumes the same  notion of  {temporal 
priority} induced by a tree. We also remark that this is also the temporal priority 
made explicit in the coefficients of Desper's method. 
Given these results, we define the following notion of causation.
\begin{definition}\label{def:our-cause} We state that $a$ is a \emph{prima facie cause} of $b$ if 
$a$ is a probability raiser of $b$,  and it occurs more frequently:
$\Pcond{b}{a} > \Pcond{b}{\~ a}$ and $\Probab{a}>\Probab{b}$.
\end{definition}

\New{We term  {\em prima facie topology} a directed acyclic graph (over
  some events) where each edge represents a prima facie cause. When at most a
  single incoming edge is assigned to each event (i.e., an event has
  at most a {\em unique cause}, in the real world), we term this
  structure {\em single-cause prima facie topology}. Intuitively, this
  last class of topologies correspond to the trees or, more generally forests when they have disconnected components, that we aim at reconstructing. 
}

%

Before moving on to introducing our algorithm let us  discuss our
definition of \emph{\cause}, \New{its role in the definition of the reconstruction problem} and some of its limitations.  As already mentioned, it may be that for some prima facie cause $c$ of an event $e$, there is a third event $a$ prior 
to both, such that $a$ causes $c$ and ultimately $c$ causes $e$. Alternatively, 
$a$ may cause both $c$ and $e$ independently, and the \cause relationship observed 
from $c$ to $e$ is merely \emph{spurious}. \New{In the context of the tree-reconstruction problem, namely when it is assumed that each event has at most a {unique cause}, the aim is to filter out the spurious edges from a general prima facie topology, so to extract a single-cause prima facie structure  (see Figure \ref{fig:prima-facie}). }


\New{Definition \ref{prop:PRtemporalPrior} summarizes Suppes basic notion of prima facie
cause, while it is ignoring deeper discussions of \cause 
that aim at distinguishing between actual genuine  and spurious
causes, e.g. screening-off, background context, d-separation \cite{cartwright1979,pearl2000,sep-causation-probabilistic}.  For our purposes however, the above  definition is sufficient when $(i)$ all the significant events are considered, i.e., all the genuine causes are observed as in  a 
closed-world assumption, and $(ii)$ we aim at extracting the \emph{order} of progression among them (or determine that there is no apparent relation), rather than  extracting causalities \emph{per se}.} Note that these assumptions  are strong and might be weakened in the future (see Discussions), but are shared by us and 
\cite{desper1999inferring}. 


%

\New{Finally, we recall a few algebraic requirements necessary for our framework to be well-defined.} First of all, the \PR{} must be computable: every 
mutation $a$ should be observed with probability strictly $0<\Probab{a}<1$. 
Moreover, we need each pair of mutations $(a,b)$ to be \emph{distinguishable} in terms of \PR{}, that 
is, for each pair of mutations $a$ and $b$, $\Pcond{a}{b}<1$ or $\Pcond{b}{a}<1$ similarly to the above condition. Any non-distinguishable pair of events can be 
merged as a single composite event.  From now on, we will assume these conditions to be verified.

\subsection*{Performance measure and synthetic datasets} \label{sec:performance}

We made use of {\em synthetic data} to evaluate the
performance of 
\New{CAPRESE} as a function of dataset size and the false positive and negative rates. Many distinct synthetic datasets  were created
for this purpose, as explained below. The algorithm's
performance was measured in terms of {\em  Tree Edit Distance} (TED, \cite{zhang1989simple}), i.e., the 
minimum-cost sequence of node edit operations (relabeling, deletion and insertion) that 
transforms the reconstructed trees into the ones generating the data. 
The choice of this measure of
evaluation is motivated by the fact that we are interested in the
\emph{structure} behind the progressive phenomenon of cancer evolution and,
in particular, we are interested in a measure of the genuine causes that we
miss and of the spurious causes that we fail to recognize (and eliminate).  
Also, since topologies with similar distributions can be structurally
different we choose to measure performance using structural distance
rather than a distance in terms of distributions.
Within the realm of `structural metrics' however, we have also evaluated the performance with the {\em Hamming Distance} \cite{hamming_distance},  another commonly used structural  metric, and we obtained analogous results (not shown here). 


\paragraph*{Synthetic data generation and experimental setting.} \label{synth} 
Synthetic datasets were generated by sampling from various random trees  constrained to have depth $\log(|G|)$, since wide branches are harder to reconstruct than straight paths, \New{and by sampling event probabilities in $[0.05,0.95]$ (see Supplementary Materials)}. 

Unless explicitly specified, in all the experiments we used $100$ distinct random trees 
(or forests, accordingly to the test to perform) of $20$ events each. This  seems a fairly reasonable number of events and is in line with the usual size of reconstructed trees, e.g. 
\cite{Dataset1,Dataset2,Dataset3,Dataset4}. 
The {\em scalability} of the techniques  was tested against the number of samples by ranging $|G|$  from $50$ to  $250$, with a step of $50$, and by replicating $10$ independent datasets for each parameters setting (see the caption of the figures for details). 

We included a form of {\em noise} in generating the datasets, in order to account for  $(i)$ the realistic presence of  {\em biological noise} (such as  the one provided by  bystander mutations, genetic heterogeneity, etc.) and $(ii)$ {\em experimental  errors}. A {noise parameter} $0 \leq \nu < 1$  denotes the probability that any event assumes a random  value (with uniform probability), after sampling from the tree-induced distribution. Algorithmically this process generates, on average, ${|G|\nu}/{2}$ random entries in each sample (e.g. with $\nu = 0.1$ we have, on average, one error per sample). 
We wish to assess whether these noisy samples can mislead the
reconstruction process, even for low values of $\nu$. Notice that
assuming a uniformly distributed noise    may appear simplistic since
some events may be more robust, or easy to measure, than
others. \New{However, introducing in the data both {\em false
    positives} (at rate $\epsilon_+=\nu/2$) and {\em negatives} (at
  rate $\epsilon_-=\nu/2$) makes the inference problem substantially
  harder, and was first investigated in \cite{szabo2002}.}

In  the Results section, we refer to datasets generated with rate $\nu >0$ as noisy synthetic dataset.  In the numerical experiments,  $\nu$ is usually discretized by  $0.025$, (i.e., $2.5\%$ noise).

\section*{Results}

\subsection*{Extracting progression trees with probability raising and a shrinkage-like estimator} \label{sec:progrextr}


\New{The CAPRESE} reconstruction method is described in Algorithm
\ref{alg:tree}. The algorithm is similar to Desper's \New{and Szabo's }algorithm, 
the main difference being an alternative weight function based on
a shrinkage-like estimator. 


\begin{definition}[Shrinkage-like estimator]\label{def:estimator} We define the {\em shrinkage-like estimator} $m_{a \to b}$ of 
the confidence in the \cause relationship from $a$ to $b$ as 
\begin{equation}
m_{a \to b} = (1-\lambda) \alpha_{a\to b} + \lambda  \beta_{a\to b}\, ,
\end{equation}
where $0\leq \lambda\leq 1$ and 
\begin{align}
\alpha_{a \to b} = \dfrac{\Pcond{b}{a} - \Pcond{b}{\~ a}}{\Pcond{b}{a} + \Pcond{b}{\~ a} }&&
\beta_{a \to b} =  \dfrac{\Pconj{a}{b} - \Probab{a}\Probab{b}}{\Pconj{a}{b} + \Probab{a}\Probab{b}}\, .
\end{align}
\end{definition}
This estimator is similar in spirit to a shrinkage estimator (see \cite{EmpiricalBayes}) and 
combines a normalized version of \PR{}, the  {\em raw estimate}
$\alpha$, with a {\em correction factor} $\beta$ (in our case a correlation-based measure of temporal
distance among events), to define a
proper order in the confidence of each \cause relationship.
Our $\lambda$ is the analogous of the {\em shrinkage coefficient} and
can have a Bayesian interpretation based on the strength of our belief that $a$ and $b$ are causally relevant to one another and the evidence that $a$ raises the probability of $b$. In the absence of a closed
form solution for the optimal value of $\lambda$, one may rely on cross-validation of simulated data. 
The power of shrinkage (and our shrinkage-like estimator) lies in the possibility of determining an optimal value for $\lambda$  to balance the effect of the correction factor on the raw model estimate \New{to ensure optimal performances on ill posed instances of the inference problem}. A crucial difference, however, between our
 estimator and classical shrinkage, is that 
our estimator aims at improving the performance of the {\em overall} reconstruction process, not limited to the  performance of the estimator 
itself as is the case in shrinkage. 
That is, the metric $m$ induces an ordering to the events reflecting
our confidence for their \cause.  
\New{Furthermore, since we make no assumption about the underlying distribution, we learn it empirically by cross-validation.}
In the next sections we show that  the shrinkage-like estimator is 
an effective  way to get such an ordering especially when \NewLast{data are} noisy. In 
\New{CAPRESE}  we use a pairwise matrix  version of the estimator.

%

\paragraph{The raw estimator and the correction factor.}

By considering only  the raw estimator $\alpha$, we would include an edge $(a \to b)$  in the tree consistently in terms of $(i)$ Definition \ref{def:our-cause}  (Methods) and $(ii)$ if $a$ 
is the best probability raiser for $b$. When $\Probab{a}=\Probab{b}$ the events $a$ and $b$ are indistinguishable  in terms of temporal priority, thus $\alpha$ is not sufficient to decide their causal relation, if any. This intrinsic ambiguity is unlikely in practice even if, in principle, it is possible. Notice that this formulation of $\alpha$  is a monotonic normalized version of the \PR{} ratio.
 \begin{prop}[Monotonic normalization] \label{prop:monotonicNorm}
For any two events $a$ and $b$ we have 
\begin{equation} 
\Probab{a} > \Probab{b} \iff
\frac{\Pcond{b}{a}}{\Pcond{b}{\~ a}} >
\frac{\Pcond{a}{b}}{\Pcond{a}{\~ b}} \iff
\alpha_{a \to b} > \alpha_{b \to a}\,.
\end{equation} 
\end{prop}
This raw model estimator satisfies  $-1\leq \alpha_{a \to b}, \alpha_{b \to a} \leq1$: when it tends to $-1$  the pair of events appear disjointly (i.e., they show an anti-\cause pattern), when it tends 
to  $0$ no \cause or anti-\cause can be inferred and the two events are statistically independent and, when it tends to $1$, the \cause relationship 
between the two events is 
\New{genuine}. 
Therefore, $\alpha$ provides a quantification of the degree of
confidence for a \PR{} \cause relationship. In fact, for any given possible \cause edge $(a,b)$, 
the term $\Pcond{b}{\~ a}$ gives an estimate of the \emph{error rate} of $b$, therefore the numerator of 
the raw model $\alpha$ provides an estimate of how often $b$ is actually caused by $a$. The $\alpha$ 
estimator is then normalized to range between $-1$ and $+1$.

However,  $\alpha$  does not provide a general criterion to disambiguate 
 \New{among genuine causes of  a given event}. We show a specific case in which $\alpha$ is not a sufficient estimator. Let us consider, for instance, a  causal linear path:  $a \to b \to c$. In this case, when evaluating the candidate parents $a$ and $b$ for $c$  we have: $\alpha_{a \to c} = \alpha_{b \to c} = 1$, \New{so $a$ and $b$ are genuine causes of  $c$, though we would like to select $b$, instead of $a$}. Accordingly, we can only infer that $t_{a}<t_{c}$ and $t_{b}<t_{c}$, i.e., a partial ordering, which does not help to disentangle the relation among $a$ and $b$ with respect to $c$. 

In this case,  the $\beta$ coefficient can be used  to determine which 
\New{of the two genuine causes occurs closer, in time, to $c$ ($b$, in the example above)}. In general, such a correction factor  provides information on the {\em temporal distance} between events, in terms of statistical dependency. In other words, the higher the $\beta$ coefficient, the closer two events are. 
Therefore, when dealing with noisy
data and limited sample sizes, there are 
situations where, by using the $\alpha$ estimator alone, we could infer a 
\New{wrong transitive edge} to be the most likely cause even in the presence of the 
\New{real} cause. For this reason, we reduce the $\alpha$ estimator to the correction factor $\beta$, which, for each given edge $(a,b)$, is 
normalized within $-1$ and $(1-\max[\Probab{a},\Probab{b}])/(1+\max[\Probab{a},\Probab{b}]) <+1$.

The shrinkage-like estimator $m$ then results in the combination of the raw \PR{} estimator $\alpha$ and of the $\beta$ correction factor, which respects the  temporal priority induced by  $\alpha$.
\begin{prop}[Coherence in dependency and temporal priority]\label{prop:betaDependency}
The $\beta$ correction factor is {\em symmetrical} and 
subsumes the same notion of dependency of the raw estimator $\alpha$, that is
\begin{equation}
\Pconj{a}{b}>\Probab{a}\Probab{b}\ \Leftrightarrow \alpha_{a \to b}>0 \Leftrightarrow \beta_{a \to b}>0
\quad \text{ and }\quad
\beta_{a \to b}=\beta_{b \to a}\,.
\end{equation}
\end{prop}


\paragraph{The independent progressions filter.}

As in Desper's approach, we also add a {\em root} $\diamond$ with $\Probab{\diamond} = 1$ in order 
to separate  different progression paths, i.e., the different sub-trees rooted in $\diamond$. 
\New{CAPRESE} initially builds a unique tree  by using  the estimator; typically, the most likely event will be at the top of the progression even if there may be rare cases where more than one event has no valid parent, in these cases we would already be reconstructing a forest. In 
the reconstructed 
tree, all the edges represent the most confident prima facie cause, although some of those could still be 
spurious causes. Then the correlation-like weight between any node $j$ and $\diamond$ is computed as
\[
\dfrac{\Probab{\diamond}}{\Probab{\diamond}+\Probab{j}}\dfrac{\Pconj{\diamond}{j}}{\Probab{\diamond}\Probab{j}} = \dfrac{1}{1+\Probab{j}} \, .
\]
If this quantity is greater than the weight of $j$ with each upstream connected element $i$, 
we consider the best prima facie cause of $j$ to be a spurious cause and we 
substitute the edge $(i \to j)$ with the edge $(\diamond \to j)$. Note
that in this work we are ignoring deeper discussions of probabilistic \cause that aim at distinguishing between actual genuine causes and spurious causes. Instead, we remove spurious causes by using a filter based on correlation because 
the probability raising of the omnipresent event $\diamond$ is not
well defined (see Methods). 
\New{In addition, we remark that} the evaluation for an edge to be a genuine or a
spurious cause takes into account all the given events. Because of
this, if events are added or removed from the dataset, the same edge
can be defined to be genuine or 
spurious 
\New{as the set of events included in the model is varied arbitrarily. }
However, since we do not consider the problem of selecting the set of
progression events, we assume that all and only the relevant events
for the problem at hand are already
\New{known a priori and included in the model.}

Finally, note that this filter is indeed implying a non-negative threshold for the shrinkage-like estimator, when an 
edge is valid. 
\begin{myth}[Independent progressions]\label{th:tree-independent}
\New{Let $G^\ast=\{a_1,\ldots,a_k\} \subset G$ a set of $k$ {\em prima facie causes} for some $b\not \in G^\ast$, and let $a^\ast=\max_{a_i \in G^\ast}\{m_{a_i \to b}\}$.  The reconstructed tree  by 
\New{CAPRESE} contains edge $\diamond \to b$ instead of $a^\ast \to b$ if, for all $a_i \in G^\ast$}
\begin{equation}
\New{\Probab{a_i,b} < \Probab{a_i}\Probab{b} \dfrac{1}{1+\Probab{b}} +  \dfrac{\Probab{b}^2}{1+\Probab{b}} \,.}
\end{equation}
%
%
\end{myth}

The proof of this theorem  can be found in the Supplementary Materials. What this theorem suggests is that, in principle, by examining the level of statistical dependency of each pair of events, it would be possible to determine how many trees compose the reconstructed 
forest. Furthermore, it suggests that 
\New{CAPRESE} could be defined by first processing the independent progressions filter, and then using $m$ to build the independent progression trees in the forest. 

To conclude, the algorithm reconstructs a well defined tree (or, more in general, forest).

\begin{myth}[Algorithm correctness]\label{th:tree-nocycles}
\New{CAPRESE} reconstructs a well defined tree $\cal T$ without disconnected components, transitive connections and cycles.
\end{myth}

\New{\NewLast{Additionally}, asymptotically with the number of samples, the reconstructed tree is the correct one.}

\begin{myth} [Asymptotic convergence]\label{th:tree-correct}
\New{Let $T=(G\cup\{\diamond\}, E, \diamond)$ be the forest  to reconstruct from a set of $s$ input samples, given as the input matrix $D$. If $D$ is strictly sampled from the distribution induced by $T$ and infinite samples are available, i.e., $s \to \infty$, 
\New{CAPRESE}  with $\lambda \to 0$ correctly reconstructs $T$.}
\end{myth} 

\New{The proof of these Theorems are also in the Supplementary
  Materials. These theorems considered datasets where the observed and
  theoretical probabilities  match, because of $s \to
  \infty$. However, data often contains false positives and negatives
  (i.e., \NewLast{data are} noisy) and resistance to their effects is desirable in
  any  inferential technique. With  this in mind, we prove a corollary of the theorem analoguos to a result appearing in \cite{szabo2002}.}


\begin{myco}[Uniform noise] \label{cor:noiseunif}
\New{Let the input matrix $D$ be strictly sampled from the
  distribution induced by $T$ \NewLast{with sample size $s \to
  \infty$}, but let it be corrupted by noise levels
  of false positives $\epsilon_+$ and false negatives $\epsilon_-$. Let $p_{\min}=\min_{x\in G}\{\Probab{x}\}$,
 \New{CAPRESE}  correctly reconstructs $T$ for $\lambda \to 0$ whenever 
\[
\epsilon_+ < \sqrt{p_{\min}}(1-\epsilon_+ - \epsilon_-)
\]
and $\epsilon_+ + \epsilon_- <1$.
}
\end{myco}

\New{
Essentially, this corollary states that 
\New{CAPRESE} (and so the estimator $m$) is robust against a  noise affecting all samples equally. Also, the fact that it holds for $\lambda \to 0$ is sound with the theory of shrinkage estimators  for which, asymptotically, the corrector factor is not needed to regularize the ill posed problem.
}

\subsection*{Optimal shrinkage-like coefficient} \label{optimal}

\New{Theorem \ref{th:tree-correct} and Corollary \ref{cor:noiseunif} state that  with infinite samples and mild constraints on the false positive/negative rates we get optimal results with  $\lambda \to 0$. Precisely, for the uniform noise model  that we applied to synthetic data (see Methods) we have  $\epsilon_+=\epsilon_-=\nu/2$, thus the hypothesis required by Corollary \ref{cor:noiseunif} is
\[
\nu < \dfrac{\sqrt{p_{\min}}}{1/2 + \sqrt{p_{\min}}}\, .
\]
For $p_{\min}=0.05$, which we set in data generation (see
Supplementary Materials), this inequality implies correct reconstruction for $\nu
< 0.3$ (a $15\%$ error rate), with infinite samples. 
However, we are interested  in performance and the optimal value of
$\lambda$ in situations in which we have finite sample sizes as well. 
Here, we empirically estimate the optimal $\lambda$ value, both in the
case of trees and forests, as a function of noise and sample size. In
the next section, we assess performance of our algorithm empirically.
}


In Figure \ref{fig:lambda}, we show the variation of the performance
of  
\New{CAPRESE} as a function of $\lambda$, for datasets
with $150$ samples generated from tree topologies. The optimal value,
i.e., lowest Tree Edit Distance (TED, see Methods), for noise-free
datasets (i.e., $\nu = 0$) is obtained for $\lambda\to0$ 
, whereas for  the  noisy datasets  a series of U-shaped curves
suggests  a unique optimum value for $\lambda\to 1/2$, immediately
observable for  $\nu < 0.15$.  Identical results are obtained when
dealing with forests (not shown here). In addition, further experiments with $n$ varying around the typical sample size ($n=150$) show that the optimal $\lambda$ is largely insensitive to the dataset size (see Figure \ref{fig:lambdaExp}). Thus we have limited our analysis to datasets with the typical sample size that is characteristic of data currently available.  



\New{Summarizing, Figures \ref{fig:lambda} and \ref{fig:lambdaExp} suggest that for sample size below $250$ without  false positives and negatives the \PR{} raw estimate $\alpha$ suffices to drive reconstruction to good results (TED is $0$ with $250$ samples), i.e.,
\begin{equation}\label{eq:shr0}
m_{a\to b} \stackrel{\lambda\to 0}{\approx} \alpha_{a\to b}  
\end{equation}
which is obtained by setting $\lambda$ to a very small value, e.g. $10^{-2}$, in order  to consider at least 
a small contribution of the correction factor too. }Conversely, when  $\nu>0$, the best performance is obtained by averaging the shrinkage-like effect, i.e.,
\begin{equation}\label{eq:shr05}
m_{a \to b} \stackrel{\lambda=1/2}{=}  \;\dfrac{\alpha_{a \to b}}{2} +  \dfrac{\beta_{a \to b}}{2}  
\,.
\end{equation}
These results suggest that, in general, a unique optimal value for the shrinkage-like coefficient can be determined, \New{even in situations not captured by   Theorem \ref{th:tree-correct} and Corollary \ref{cor:noiseunif}.}

\subsection*{Performance of \New{CAPRESE} compared to {\em oncotrees}} \label{comparison}

\New{An analogue of  Theorem \ref{th:tree-correct}  holds for Despers's oncotrees (Theorem 3.3, \cite{desper1999inferring}), and an analogue  of Corollary \ref{cor:noiseunif}
holds for Szabo's extension with uniform noise (Reconstruction Theorem
1, \cite{szabo2002}). 
Thus, with infinite samples both approaches reconstruct the correct
trees/forests. 
With
finite samples and noise, however,
their performance may show different patterns, as speed of convergence
may vary. 
We investigate this issue in the current section. 
}

In Figure \ref{fig:treesynthnfree}  we compare the performance of
\New{CAPRESE} with  oncotrees, for the case of noise-free
synthetic data with the optimal shrinkage-like coefficient: $\lambda \to 0$, equation (\ref{eq:shr0}). \New{Since Szabo's algorithm is equivalent to Desper's
  without false negatives and positives, we rely solely on Szabo's
  implementation \cite{szabo2002}.}
In Figure \ref{fig:exampleTree} we show an example of reconstructed tree where 
\New{CAPRESE} infers the correct tree while oncotrees mislead a \cause relation.

In general, one can observe that the TED of 
\New{CAPRESE}
is, on average, always bounded above by the TED of oncotrees, both in
the case of trees and forests. For trees,  with $50$ samples  the
average TED of 
\New{CAPRESE} is around $6$, whereas for
Desper's technique it is around $13$. The performance of both
algorithms improves as long as the number of samples is increased:
\New{CAPRESE} has the best performance (i.e., TED $\approx
0$)  with $250$ samples, while  oncotrees have TED around $6$. When
forests are considered, the difference between the performance of the
algorithms reduces slightly, but also in this case 
\New{CAPRESE} clearly {outperforms oncotrees}. 

Notice that the improvement due to the increase in the sample size seems to reach a {\em plateau},  and the initial TED for our  estimator seems rather close to the plateau value. 
\New{This empirical analysis suggests that  
\New{CAPRESE} has already good performances with few samples, a favorable adjoint to Theorem \ref{th:tree-correct}.}  This result has some important practical implications, particularly considering  the scarcity of available biological data.

In Figure \ref{fig:synteticnoisy}  we extend the comparison to  {\em
  noisy} datasets. In this case, we used the optimal shrinkage-like
coefficient: $\lambda \to 1/2$, equation (\ref{eq:shr05}). The results
confirm what observed without false positives and negatives, as
\New{CAPRESE} outperforms oncotrees up to $\nu = 0.15$, for
all the sizes of the sample sets.
In the Supplementary Materials we show similar plots for the noise-free case.

\New{We can thus draw the conclusion that our algorithm performs better  with finite samples and noise, since less samples are required to get good performances and a higher resistance to false positives and negatives is shown.}

\New{
\subsubsection*{Performance of CAPRESE compared to
  \emph{Conjunctive Bayesian Networks}}
Inspired by Desper's seminal work, Beerenwinkel and others developed
methods to  estimate the constraints on the order in which mutations
accumulate during cancer progression, using a probabilistic graphical model
called  \emph{Conjuntive Bayesian Networks} (CBN)\cite{beerenwinkel2007conjunctive,gerstung2009quantifying}. 
While the goal of this research was to reconstruct \emph{direct acyclic
  graphs} and not trees per se, evidence presented in
\cite{Hainke:2012hwa} suggests that, in the absence of noise, these models perform better than
oncotrees even at reconstructing \emph{trees}. For this reason, we performed
experiments similar to the ones suggested above, comparing 
\New{CAPRESE}
to the extension of CBN called \emph{hidden-CBN} (h-CBN) that accounts for noisy
genotype observations \cite{gerstung2009quantifying}. 
This method combines CBNs with a simulated
annealing algorithm for structure search and a denoising of the
genotypes via the maximum a posteriori estimates to compute the
most likely progression. \NewLast{One aspect that complicates a
  comparison between CAPRESE and (h-)CBN is that the methods assume 
  different models. For example, at the heart of CBN is a  monotonicity assumption
  (i.e., an event can only occur if all its predecessors have
  occurred) not assumed by CAPRESE. Despite the differences between
  the model assumptions, we present a
preliminary comparison between the methods in Appendix
\ref{sec:figures}, indicating that we not only outperform
oncotrees, but h-CBNs as well.} In particular, this suggests that
\New{CAPRESE} converges much faster than h-CBNs with
respect to the sample size, also in the presence of noise.

\NewLast{We also analyze the rate of
  \emph{false positives/negatives} reconstructed by CAPRESE  when (synthetic)
  data are sampled from DAGs (Appendix
  \ref{sec:figures}). The rate of \emph{false positives} goes
  to $0$ as the sample size increases, implying that CAPRESE is capable
  of reconstructing a tree subsumed by the underlying causal DAG
  topology. In addition, the number of \emph{false negatives}
  approaches a value proportional to the connectivity of the model from
  which the data was generated. This is expected, since CAPRESE will  assign at most one cause to each 
considered event. } 
However, it should be noted that further
experiments \NewLast{with samples selected from a wider array of topologies} should be performed to
confirm these results and compare both methods in full. While not within the scope of the
current paper, these issues will be
addressed in future work.
}

\subsection*{Case studies}\label{sec:realdata}

In the next subsections we apply 
\New{CAPRESE} to real cancer data obtained by    {\em Comparative Genomic Hybridization}  (CGH) and {\em Next Generation Sequencing}  (NGS). This shows the potential application of reconstruction techniques to various types of mutational profiles and various cancers.

\subsubsection*{Performance on cancer CGH datasets}  

\New{Encouraged by the
results in previous sections}, we test our reconstruction approach on a real {\em ovarian cancer} dataset made available within the oncotree
package \cite{desper1999inferring}. The data was collected through the public 
platform SKY/M-FISH \cite{knutsen2005interactive}, used to allow investigators to  share  molecular 
cytogenetic data. The data was obtained by using the CGH technique on samples from {\em papillary serous cystadenocarcinoma} 
of the ovary. 
This technique uses fluorescent staining to detect CNV data at the resolution of
chromosome arms. 
While the recent emergence of NGS approaches make
the dataset itself rather outdated, 
the underlying principles remain the same and the dataset provides a
valid test-case for our approach.
The seven most
commonly occurring events are selected from the $87$ samples, and the
set of events are the following gains and losses on chromosomes arms
$G=\{ 8q+, 3q+, 1q+, 5q-, 4q-, 8p-, Xp-\}$ (e.g.,  $4q-$ denotes
a deletion of the $q$ arm of the $4^{th}$ chromosome). 

In Figure \ref{fig:treesOvarian} we compare the trees reconstructed by the two approaches.
Our technique differs from Desper's   by predicting the causal sequence of alterations 
\[
8q+\;\to\;8p-\;\to \;Xp-\, ,
\]
\New{when used either $\lambda \to 0$ or $\lambda=1/2$. Notice that among the samples in the dataset some are not generated by the distribution induced by the recovered tree, thus comparing the reconstruction for both $\lambda$s is necessary.}


At this point, we do not have a biological interpretation for this
result. However, we do know that common cancer genes reside in these
regions,  e.g. the tumor suppressor gene \textsc{Pdgfr} on $5q$  and
the oncogene \textsc{Myc} on $8q$), and loss of heterozygosity on the
short arm of chromosome $8$ is 
\New{quite}
common (see, e.g., \url{http://www.genome.jp/kegg/}). Recently, evidence has been reported that  $8p$ contains many cooperating cancer genes  \cite{xue2012cluster}.  

In order to assign a confidence level to these inferences we applied 
both parametric  and non-parametric  {\em bootstrapping methods} to our results. Essentially, these tests consist of using the reconstructed trees (in the
parametric case), or the probability observed in the dataset (in the
non-parametric case) to generate new synthetic datasets, and then
reconstructs again the progressions (see, e.g., \cite{efron1982jackknife}
for an overview of these methods \New{and \cite{Szabo2008} for the use of bootstrap for evalutating the 
confidence of oncogenetic trees.}). 
The confidence is given by the number of times the trees in Figure
\ref{fig:treesOvarian} are reconstructed from the generated data.  
A similar approach can be used to estimate the confidence of every
edge separately. 
For oncotrees the \emph{exact tree} is obtained $83$ times
out of $1000$ non-parametric resamples, so its estimated confidence is
$8.3\%$. For 
\New{CAPRESE} the confidence is $8.6\%$. In the
parametric case \New{with false positive and false negative error
  rates of $0.21$ and $0.027$, following  \cite{szabo2002},} the
confidence of oncotrees is $17\%$ while the confidence of our method is much
higher: $32\%$. \New{When error rates are forced to $0$ the confidence
  of oncotrees raises to $86.6\%$ and $90.9\%$ respectively.}

For the non-parametric case, edges confidence is shown in Table
\ref{fig:freqOvarian}. 
Most notably, our algorithm reconstructs the inference
$8q+\to8p-$ with high confidence (confidence $62\%$, and $26\%$ for $5q-\to8p-$),
while the confidence of the edge $8q+\to5q-$ is only $39\%$, almost
the same as $8p- \to 8q+$ (confidence $40\%$).  \New{The confidences are similar with either $\lambda \to 0$ or $\lambda=1/2$}.

\paragraph{Analysis of other CGH datasets.}

We report the differences between the reconstructed trees also based on
datasets of gastrointestinal and oral cancer (\cite{Dataset2,Dataset4} respectively). In the case of gastrointestinal stromal  cancer, among the $13$ CGH events considered in \cite{Dataset2}
(gains on $5p$, $5q$ and  $8q$, losses on $14q$, $1p$, $15q$, $13q$, $21q$, $22q$, $9p$, $9q$,  $10q$ and $6q$), oncotrees identify the path progression 
\[
1p- \to 15q- \to 13 q- \to 21q-
\] 
 while 
 \New{CAPRESE} reconstructs the branch 
\begin{align*}
1p- \to 15q- && 1p- \to 13q- \to 21q-\,.
\end{align*}
In the case of  oral cancer, among the $12$ CGH events considered in \cite{Dataset4}
(gains on $8q$, $9q$, $11q$, $20q$, $17p$, $7p$, $5p$, $20p$ and  $18p$, losses on $3p$, $8p$ and $18q$), the reconstructed trees differ since oncotrees identifies the  path 
\begin{align*}
8q+ \to 20q+  \to 20p+ 
\end{align*}
 while our algorithm reconstructs the path 
\begin{align*}
3p- \to 7p+  \to 20q+ \to 20p+\,.
\end{align*}

These examples show that 
\New{CAPRESE} provides important differences in the reconstruction compared to oncotrees.

\subsection*{Performance on cancer NGS datasets} \label{sec:realdata2} 

\newcommand{\setbp}{\textsc{Setbp1}}
\newcommand{\asxl}{\textsc{Asxl1}}
\newcommand{\tet}{\textsc{Tet2}}
\newcommand{\kras}{\textsc{Kras}}
\newcommand{\aCML}{\textsc{aCML}}

In this section we show the application of reconstruction techniques to the validation of a specific relation among recurrent mutations involved  in {\em atypical Chronic Myeloid Leukemia} (\aCML{}).

In  \cite{NatGenRP} Piazza {\em et al.} used  high-throughput {\em exome sequencing technology} to identity somatically acquired mutations in $64$ \aCML{} patients, and found a previously unidentified recurring {\em missense point  mutation}  hitting \setbp{}. By re-sequencing \setbp{} in samples with \aCML{} and other common human cancers, they found that around $25\%$ of the \aCML{} patients tested positive for \setbp{}, while most of the other types of tumors were negative.  
Assessing the possible relationship between \setbp{} variants and the mutations in many driver \aCML{} oncogenes such as (e.g., \asxl{}, \tet{}, \kras, etc.) no significant association or mutual exclusion with \setbp{} was found but for \asxl{}, which was frequently mutated together with \setbp{}, hinting at a potential relation among the events. In particular, \asxl{} was presenting either a  {\em non-sense point} or a {\em indel} type of somatic mutation.

Hence, we reconstructed \aCML{} progression models from the datasets provided in \cite{NatGenRP}, with the goal of 
\New{assessing} a {\em potential causal dependency} between mutated  \setbp{} and \asxl{}. A more extensive analysis is postponed, as we only seek to clearly illustrate the functionalities of the
algorithm here. 

As a first case (Figure \ref{fig:treesACML}, left), we
\New{treated} the \asxl{} {missense
  point} and {indel} mutations \New{as indistinguishable, and} we merged the two events in the dataset.  Afterwards, we separated the two types of mutations for \asxl{} (Figure \ref{fig:treesACML}, right).

In particular, it is interesting to notice that, when the \asxl{} mutations are considered equivalent, the inference suggests that the mutations belong to two independent progression paths (i.e., the independent progression filters ``breaks'' every potential causal relation among the events).  Conversely, when the mutations are kept separate, the progression 
model suggests that: $(i)$ the missense point  mutation  hitting \setbp{} can cause a non-sense point mutation in \asxl{} and $(ii)$ the observed \asxl{} mutations seems to be independent. Concerning edges confidence, as before assessed via non-parametric bootstrap, it is worth noting that the confidence in the indel \asxl{} mutation being an early event raises consistently in the latter case.

All in all, it seems that a progression model 
\New{allows to test the significance of} the association firstly observed in \cite{NatGenRP} and also refines the knowledge by suggesting a specific causal  and  temporal relations among events. With this this in mind, ad-hoc sequencing experiments might be set up to 
\New{assess} these predictions, eventually providing a strong evidence that could be used to, e.g., synthesize  a progression-specific \aCML{}-effective drug.

\section*{Discussion}\label{sec:concl}

In this work we have introduced a novel  framework for the reconstruction of the causal topologies  underlying cumulative progressive phenomena, based on the {\em probability raising} notion of  \cause.  Besides such a  probabilistic notion, we also introduced the use of  a {\em shrinkage-like estimator} to efficiently unravel ambiguous causal relations, often present when data are noisy. 
As a first step towards the definition of our new framework, we have
here presented an effective novel technique \New{called CAPRESE (CAncer PRogression Extraction with Single Edges)} for the reconstruction of
tree or, more generally, forest models of progression which combines
probabilistic \cause and a shrinkage-like estimation. 

\New{We prove correctness of CAPRESE by showing asymptotic
  convergence to the correct tree. Under mild constraints on the noise
  rates, this result holds for the reconstruction problem in the
  presence of uniform noise as well. Moreover, we also compare our technique
to the standard tree reconstruction algorithm based on correlation
(i.e., Oncontrees) and to a more general bayesian probabilistic
graphical model (i.e., Conjunctive Bayesian Networks)}, and show that
\New{CAPRESE} outperforms the state-of-the-art on synthetic data, also exhibiting a noteworthy efficiency with relatively small datasets. Furthermore, we  tested our technique on  ovarian, gastrointestinal and oral cancer CGH data and NGS leukemia data.
The CGH analysis suggested that our approach can infer, with high
confidence, novel causal relationships which would remain
\New{unpredictable in
a correlation-based attack}. The NGS analysis allowed validating a causal and temporal relation among key mutations in atypical chronic myeloid leukemia.

One of the strong points of 
\New{CAPRESE} is that it can be applied to genomic data of any kind, even heterogeneous, and at any resolution, as shown. In fact, it simply requires a set of samples in which the presence or the absence of some alterations supposed to be involved in a causal cumulative process have been assessed. 
Notice also that the results of our technique  can be used not
only to describe the \emph{progression} of the process, but also to
\emph{classify} different progression types. In the case of cancer,
for instance, this genome-level classifier could be used to  group
patients \New{according to the position of the detected individual
  genomic alterations in the progression model (e.g., at a specific
  point of the tree)} and, consequently, to set up a genome-specific
\emph{therapy design} \New{aimed, for instance, at blocking or slowing
  certain progression paths instead of others, as was studied in \cite{olde2014cancer}}. 

Several future research directions are possible. First, more  complex models of progression, e.g. directed acyclic graphs,  could be inferred with probability raising and compared to the standard approaches of \cite{beerenwinkel2007conjunctive,gerstung2009quantifying,gerstung2011temporal}, as we explained in the introduction. These models, rather than trees, could explain  the common phenomenon of {\em preferential progression paths} in the target process via, e.g., {\em confluence} among events. In the case of cancer, for instance, these models would be certainly more suitable than trees to describe the accumulation of mutations. 

Second, the shrinkage-like estimator itself could be improved by introducing, for instance, different correction factors. In addition, an analytical formulation of the optimal shrinkage-like coefficient could be investigated by starting from the hypotheses which apply to our problem setting, along the lines of \cite{JSestimator}. 

Third, advanced statistical techniques such as \emph{bootstrapping}
\cite{efron1982jackknife} could be used to account for   more
sophisticated models of noise within  data, so as to  decipher complex  causal dependencies.  

Finally, a further development of the framework could  involve the
introduction of {\em timed data}, in order to extend  our techniques to settings where a temporal information on the samples is available. 

\New{\subsection*{Software availability}
The implementation of CAPRESE is part of the TRanslational ONCOlogy
(TRONCO) R package and is available for download at standard R repositories.
}

\New{\subsection*{Acknowledgements}
We are grateful for the many excellent comments we received from anonymous
reviewers. }

\bibliographystyle{ieeetr}
\bibliography{progression}

\newpage
\appendix


\section{Figures and tables}
\begin{figure}[!htbp]
\begin{center}
 \includegraphics[width=0.60\textwidth]{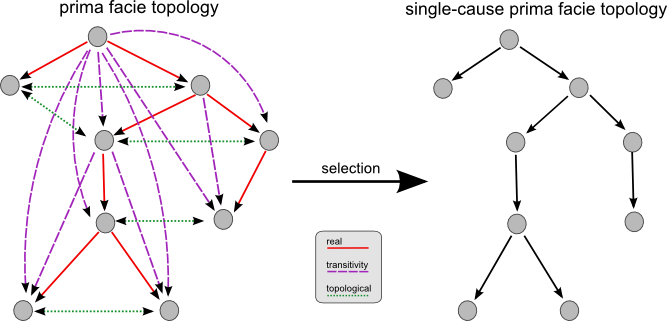}
\end{center}
\caption[Prima facie topology]{\New{{\bf Prima facie topology.}
    Example prima facie topology where all edges $(a,b)$  represent prima facie causes, according to Definition \ref{def:our-cause}: $a$ is a probability raiser of $b$  and it occurs more frequently. In left, we filter out spurious causes and select only the real ones among the genuine, yielding a single-cause prima facie topology.} }
\label{fig:prima-facie}
\end{figure}

\begin{algorithm}[!htbp]
\caption{CAPRESE: a tree-like reconstruction with a shrinkage-like estimator}
\label{alg:tree}
\begin{algorithmic}[1]
\STATE  consider a set of $n$ genetic events $G$ plus a special event $\diamond$, added to each sample of the dataset;
\STATE define a  $m  \times n$ matrix $M$ where each entry contains the shrinkage-like estimator
\[
m_{i \to j} = (1-\lambda)\cdot {\dfrac{\Pcond{j}{i} - \Pcond{j}{\~i}}{\Pcond{j}{i} +
  \Pcond{j}{\~i}} + \lambda \cdot \dfrac{\Pconj{i}{j} - \Probab{i}\Probab{j}}{\Pconj{i}{j} + \Probab{i}\Probab{j}}}
\]
according to the observed probability of the events $i$ and $j$;
\STATE [\PR{} \cause{}] define a tree ${\cal T}=(G\cup\{\diamond\},E, \diamond)$ where $(i \to j) \in E$ for $i,j\in G$ if and only if:
\[
m_{i \to j}>0
\quad \text{and} \quad
m_{i \to j} > m_{j \to i} 
\quad \text{and}\quad \forall i' \in G, \: m_{i,j} > m_{i',j} \, .
\]
\STATE [Independent progressions filter] define $G_j=\{x \in G\mid \Probab{x} > \Probab{j} \}$, replace edge $(i \to j) \in E$ with edge 
$(\diamond \to  j)$ if, for all $x \in G_j$, it holds 
\[
\dfrac{1}{1+\Probab{j}} > \dfrac{\Probab{x}}{\Probab{x}+\Probab{j}}\dfrac{\Pconj{x}{j}}{\Probab{x}\Probab{j}}
\,.\]
\end{algorithmic}
\end{algorithm}

\begin{figure}[!htbp]
\begin{center}
 \includegraphics[width=0.50\textwidth]{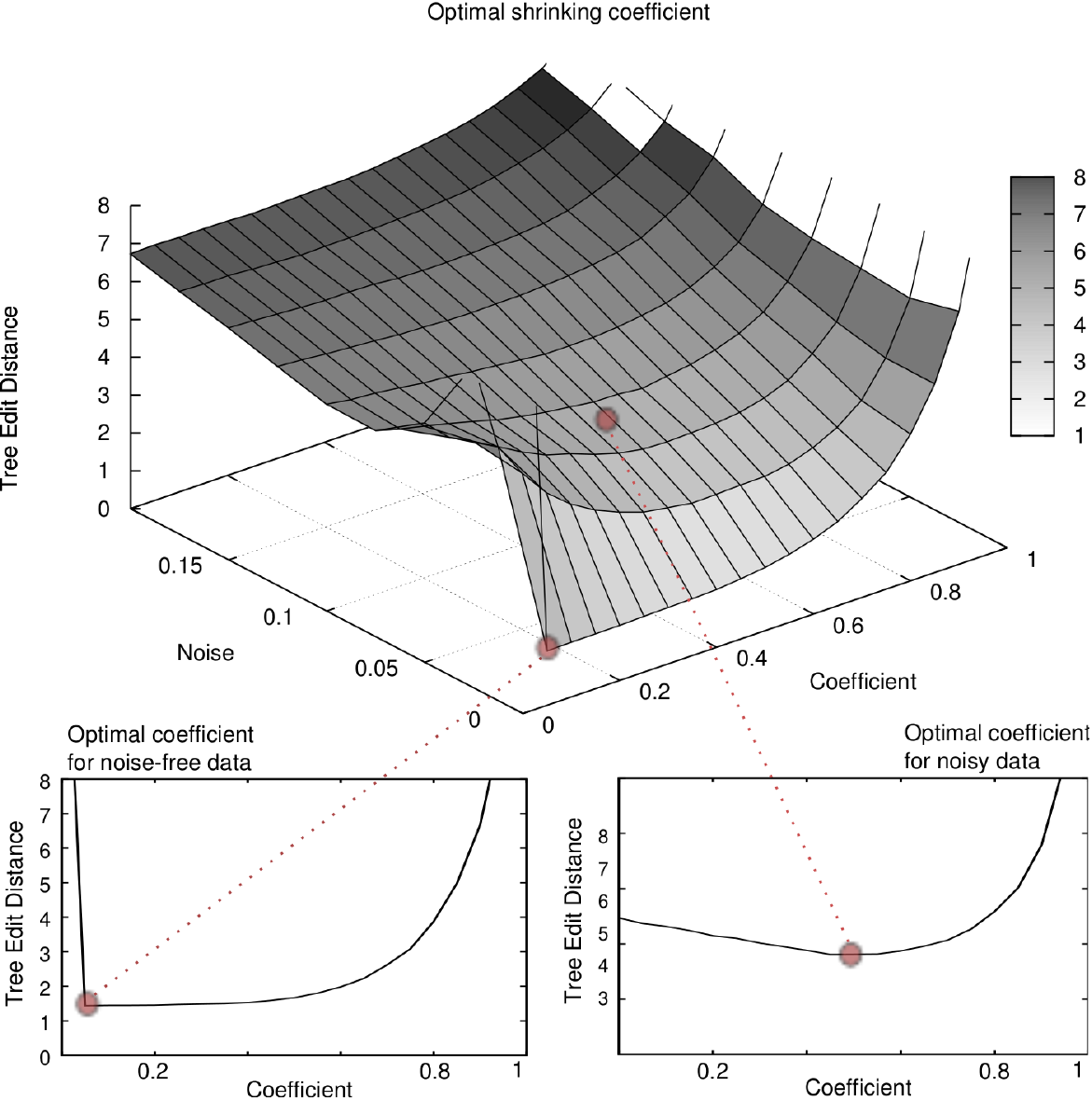}
\end{center}
\caption[Optimal shrinkage-like coefficient for reconstruction performance]{{\bf
    Optimal shrinkage-like coefficient for reconstruction performance.}
    We show here the performance in the reconstruction of trees (TED surface) with $n=150$ samples as a function of the shrinkage-like coefficient $\lambda$. Notice the global optimal performance for 
    $\lambda\to0$ when $\nu\to0$ and for $\lambda\approx1/2$ when $\nu>0$.}
\label{fig:lambda}
\end{figure}

\begin{figure}[!htbp]
\begin{center}
	\includegraphics[width=0.48\textwidth]{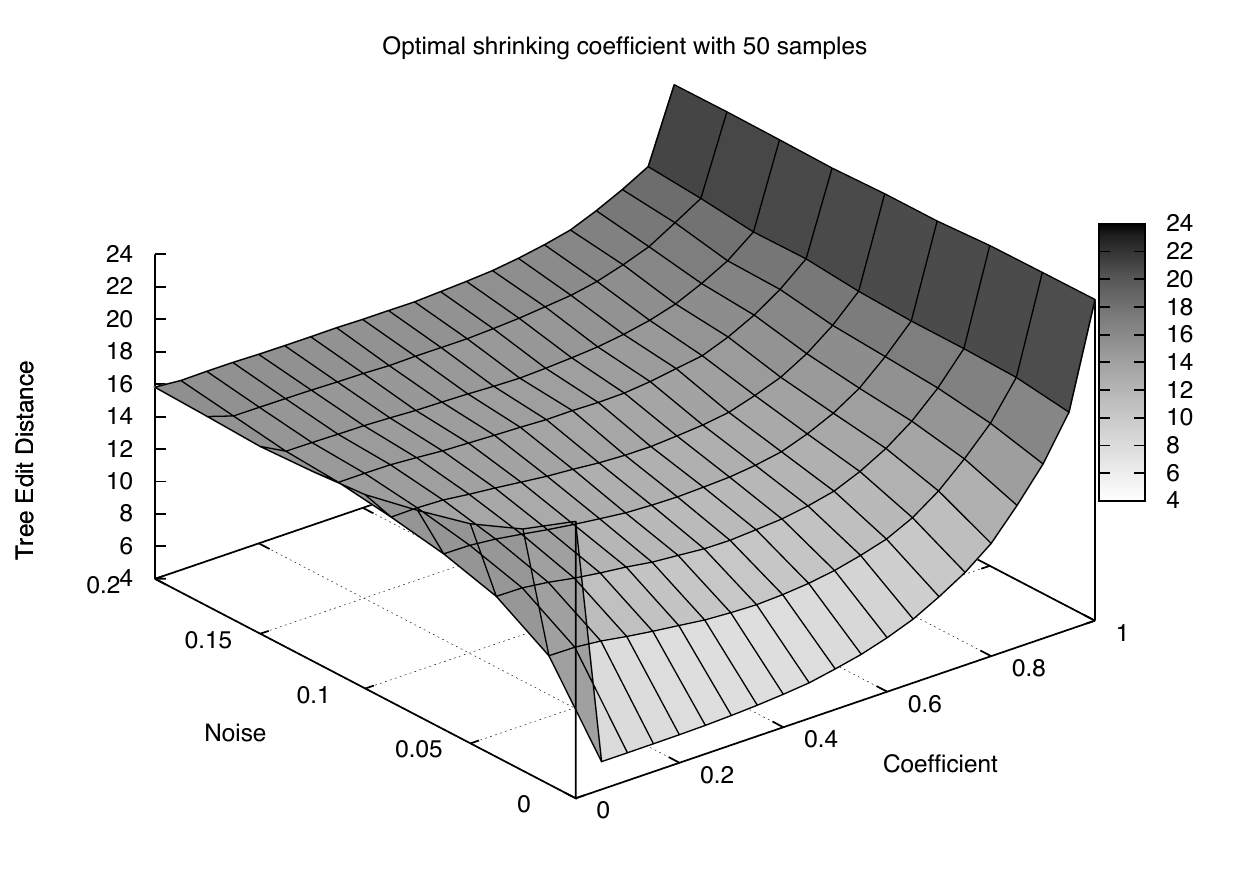} 
	\includegraphics[width=0.48\textwidth]{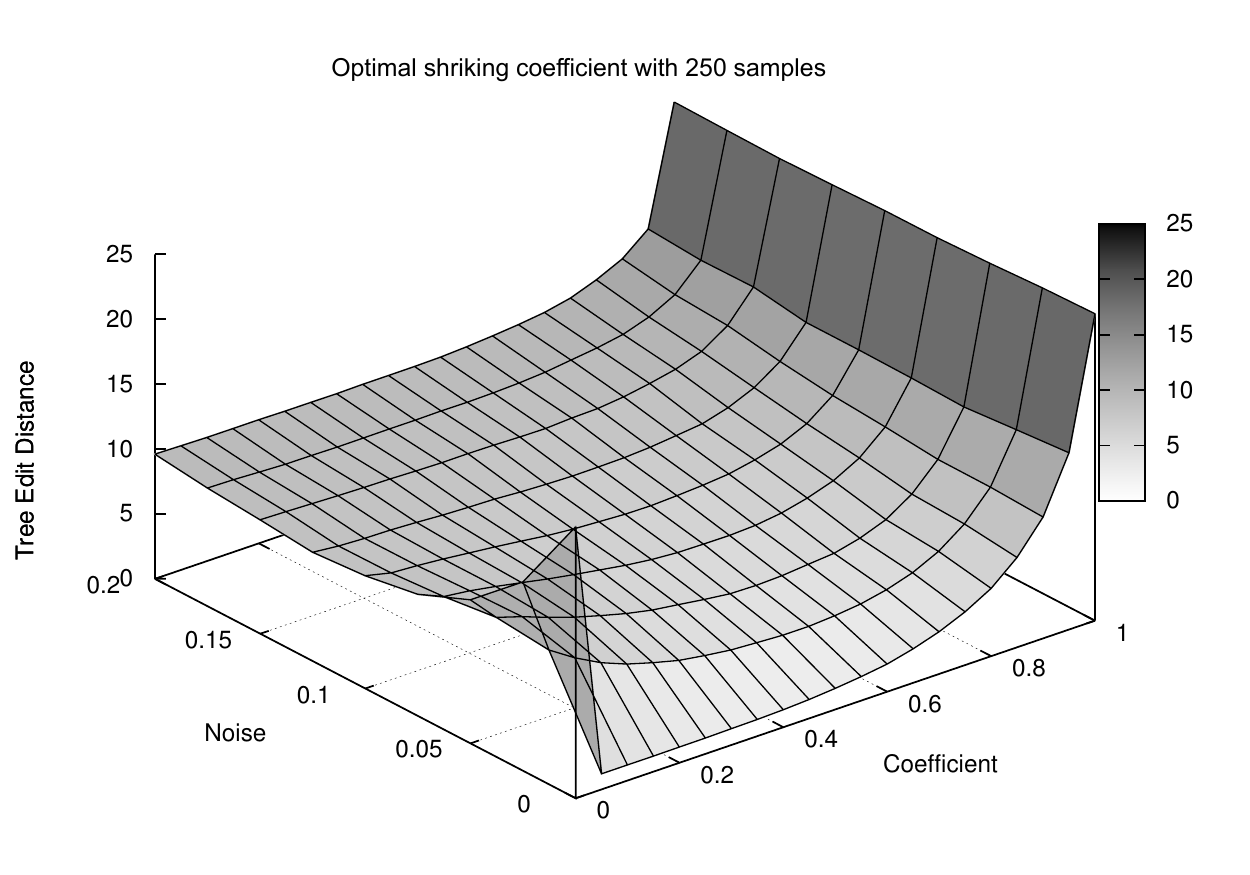}
\end{center}
\caption[Optimal $\lambda$ with datasets of different size.]{{\bf Optimal $\lambda$ with datasets of different size.}
We show the analogous of Figure \ref{fig:lambda} with $50$ and $250$ samples. The estimation of the optimal shrinkage-like coefficient $\lambda$ appears to be irrespective of the sample size.}
\label{fig:lambdaExp}
\end{figure}

\begin{figure}[!htbp]
\begin{center}
 \includegraphics[width=0.57\textwidth]{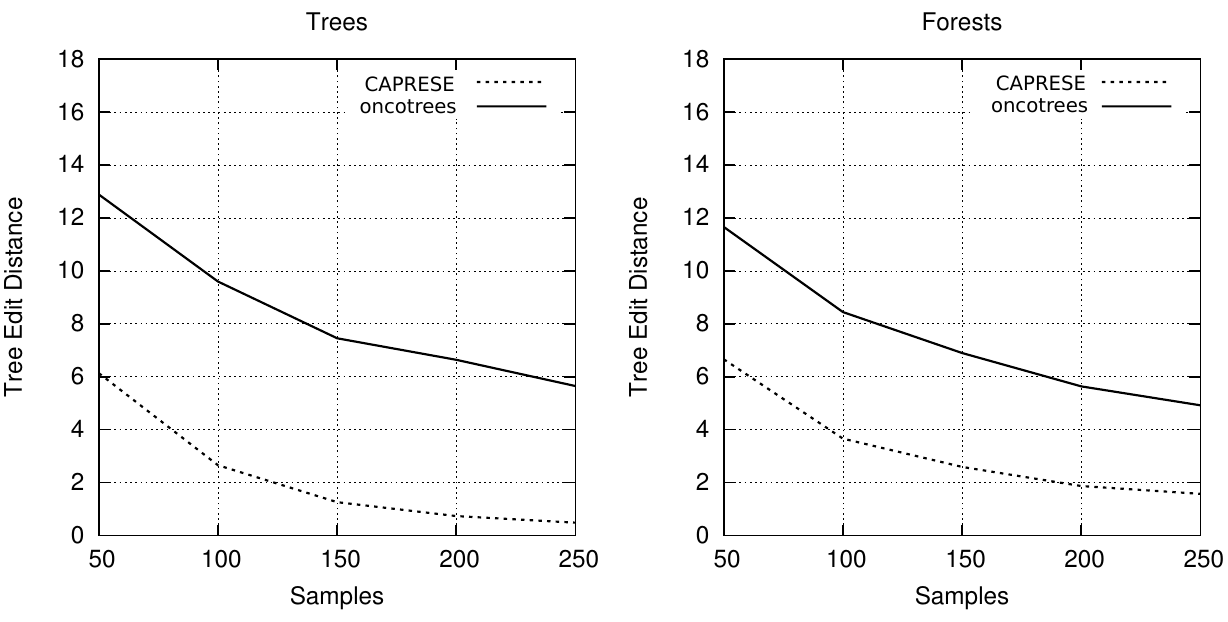}
\end{center}
\caption[Comparison on noise-free synthetic data]{{\bf
    Comparison on noise-free synthetic data.}
Performance of  
\New{CAPRESE} (dashed line) and oncotrees (full line) in average TED when data are generated by random trees (left) and forests (right). 
In this case $\nu=0$ (no false positives/negatives) and $\lambda\to 0$ in the estimator $m$.
}
\label{fig:treesynthnfree}
\end{figure}

\begin{figure}[!htbp]
\begin{center}
\includegraphics[width=0.59\textwidth]{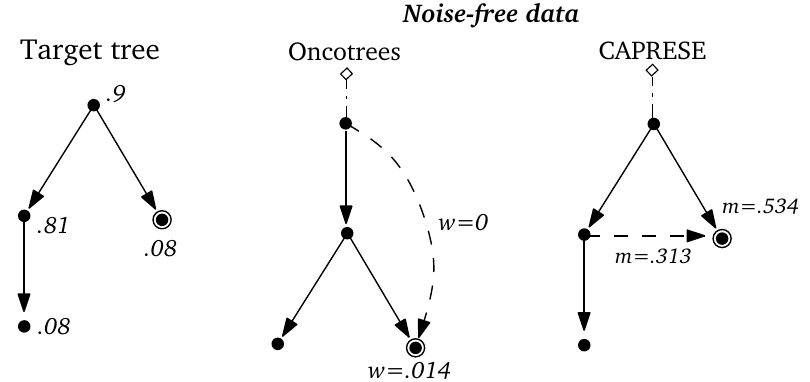}
\end{center}
\caption[Example of reconstructed trees.]{{\bf Example of reconstructed trees.}  
Example of reconstruction from a dataset with \New{$100$ samples} generated by the left tree \New{(the theoretical probabilities are shown, i.e., the doubly-circled event appears in a sample with probability $.08$)}, with $\nu=0$. \New{In the sampled dataset oncotrees mislead the cause of  the doubly-circled mutation}  ($w=0$ for the true edge and $w=0.014$ for the wrong one). 
\New{CAPRESE} infers the correct cause (the values of the estimator $m$ \New{with $\lambda=1/2$} are shown, similar results are obtained for $\lambda\to 1$).  
  }
\label{fig:exampleTree}
\end{figure}

\begin{figure}[!htbp]
\begin{center}
\includegraphics[width=0.99\textwidth]{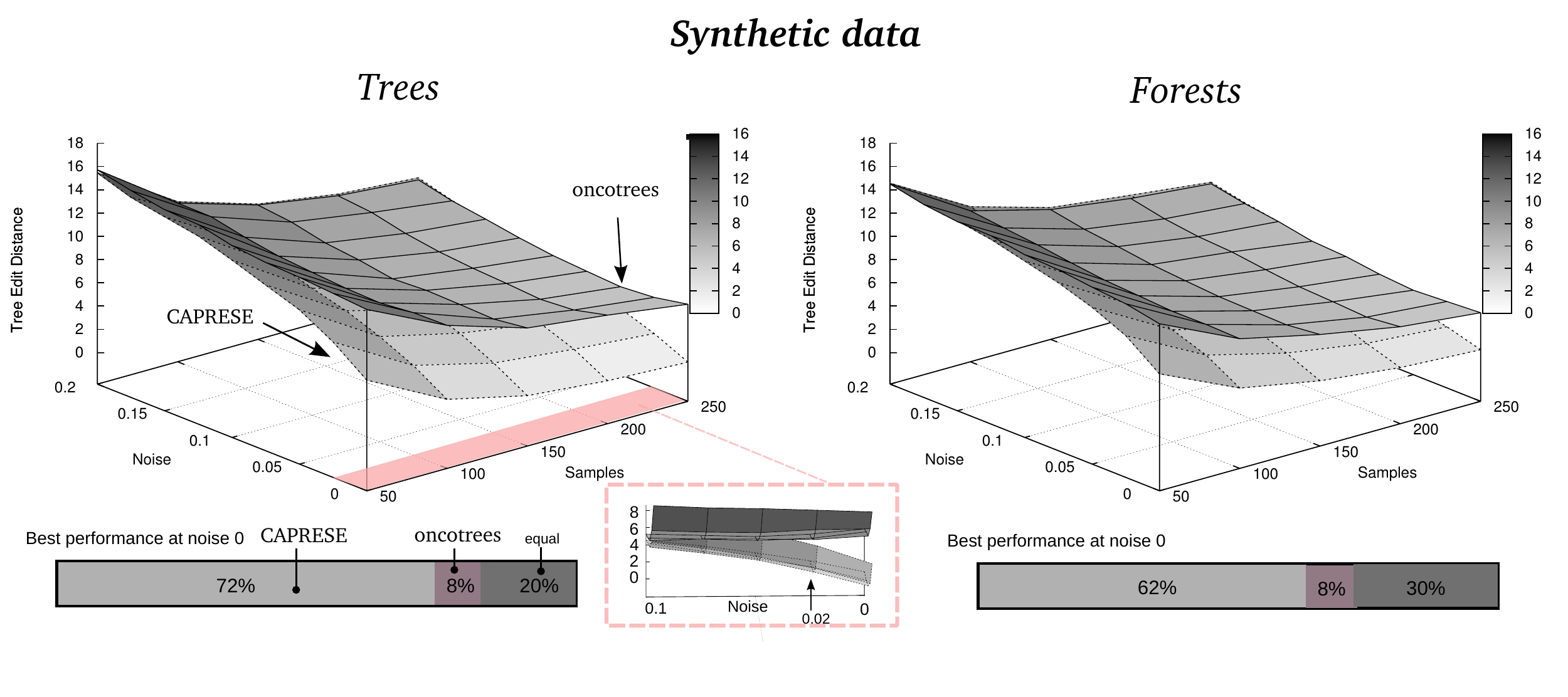}
\end{center}
\caption[Reconstruction with noisy synthetic data and $\lambda=1/2$.]{{\bf Reconstruction with noisy synthetic data and $\lambda=1/2$.} Performance of 
\New{CAPRESE} and oncotrees as a function of the number of samples and noise $\nu$. According to Figure \ref{fig:lambda} the shrinkage-like coefficient is set to $\lambda=1/2$. The magnified image shows the convergence to Desper's performance for $\nu \approx 0.1$. 
The barplot represents the percentage of times the best performance is achieved at $\nu=0$.  
  }
\label{fig:synteticnoisy}
\end{figure}

\begin{figure}[!htbp]
\begin{center}
\includegraphics[width=0.75\textwidth]{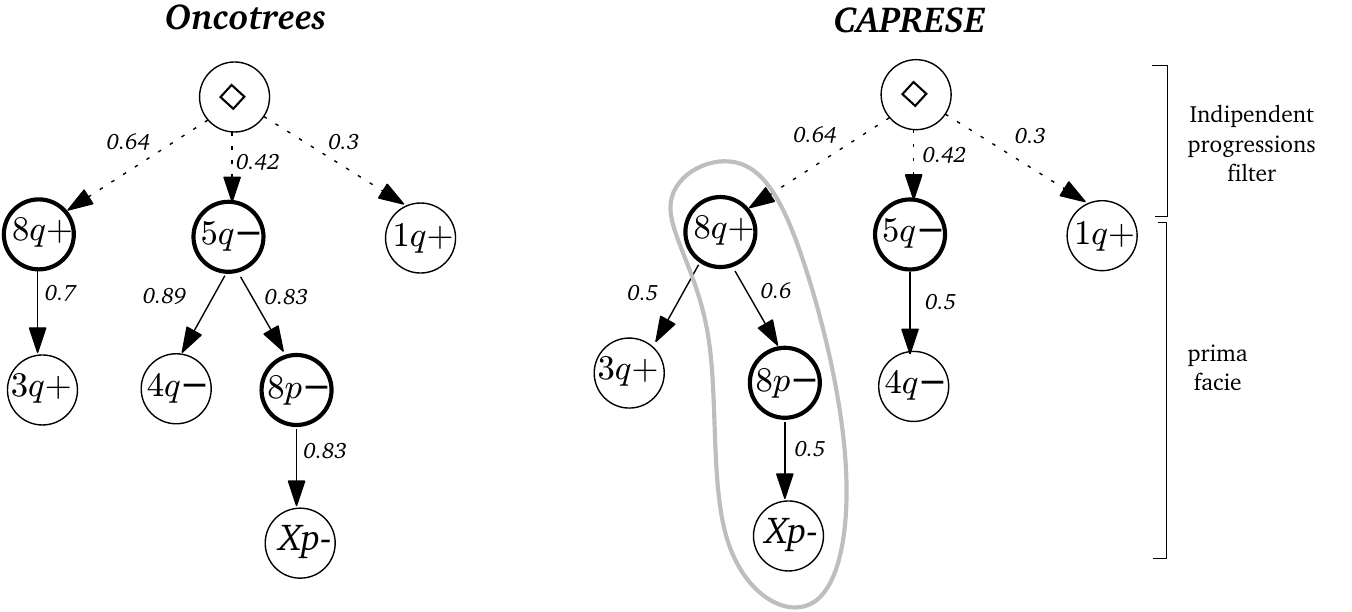}
\end{center}
\caption[Tree reconstruction of ovarian cancer progression]{{\bf Reconstruction of ovarian cancer progression.} Trees reconstructed by oncotrees and  
\New{CAPRESE} (with $\lambda \to 0$, \New{with $\lambda=1/2$ the same tree is reconstructed}). The set of CGH events considered are
gains on $8q$, $3q$ and $1q$ and losses on $5q$, $4q$, $8p$ and $Xp$. Events on chromosomes arms containing the key genes for ovarian cancer are in bolded circles. In the left tree all edge weights are the observed probabilities of events. In the right the full edges are the \cause inferred with the \PR{} and the weights represent the scores used by 
\New{CAPRESE}. Weights on dashed lines are as in the left tree.
}
\label{fig:treesOvarian}
\end{figure}

\begin{figure}[!htbp]  
\begin{center} \footnotesize Oncotrees (overall confidence $8.3\%$)\\
  \begin{tabular}{  c  | c |  c |  c |  c |  c |  c |  c | }
    $\to$  & $8q+$ & $3q+$ & $5q-$ & $4q-$ & $8p-$ & $1q+$ & $Xp-$ \\ \hline
    $\diamond$ & \hC{$\mathbf{.99}$} & \lC{$.06$} & \mC{$\mathbf{.51}$} & \lC{$.22$} & \lC{$.004$} & \hC{$\mathbf{.8}$} & \lC{$.06$} \\ \hline
     $8q+$  & $0$ & \hC{$\mathbf{.92}$} & \lC{$.08$} & \lC{$0.16$} & \mC{$0.4$} & \lC{$.02$} & \lC{$.007$}  \\ \hline
    $3q+$ &  \lC{$.002$} & $0$ & \lC{$.04$} & $0$ & $0$ & \lC{$.09$} & \lC{$.04$}  \\ \hline
        $5q-$& \lC{$.001$} & \lC{$.002$} & $0$ & \mC{$\mathbf{.52}$} & \lC{$\mathbf{.39}$} & \lC{$.009$} & \lC{$.16$}  \\ \hline
    $4q-$  &$0$ &$0$ &\lC{$.27$} &$0$ &\lC{$.14$} &\lC{$.05$} &\lC{$.11$}   \\ \hline
    $8p-$  & $0$ & $0$ & \lC{$.07$} & \lC{$.08$} & $0$ & \lC{$.004$} & \mC{$\mathbf{.59}$} \\ \hline
        $1q+$ & $0$ & $0$ & $0$ & \lC{$.004$} & $0$ & $0$ & $0$  \\ \hline
    $Xp-$ & $0$ & $0$ & $\lC{.003}$ & $\lC{.003}$ & \lC{$.04$} & \lC{$.01$} & $0$ \\ \hline
      \end{tabular}\medskip

 \New{CAPRESE} (overall confidence $8.6\%$)\\
  \begin{tabular}{  c  | c |  c |  c |  c |  c |  c |  c | }
    $\to$  & $8q+$ & $3q+$ & $5q-$ & $4q-$ & $8p-$ & $1q+$ & $Xp-$ \\ \hline
    $\diamond$ & \hC{$\mathbf{.99}$} & \lC{$.06$} & \mC{$\mathbf{.51}$} & \lC{$.22$} & \lC{$.004$} & \hC{$\mathbf{.8}$} & \lC{$.06$} \\ \hline
     $8q+$  & $0$ & \hC{$\mathbf{.92}$} & \lC{$.06$} & \lC{$.16$} & \mC{$\mathbf{.62}$} & \lC{$.01$} & \lC{$.008$}  \\ \hline
    $3q+$ &  \lC{$.002$} & $0$ & \lC{$.03$} & $.002$ & $0$ & \lC{$.09$} & \lC{$.04$}  \\ \hline
        $5q-$& \lC{$.001$} & \lC{$.002$} & $0$ & \mC{$\mathbf{.5}$} & \lC{${.26}$} & \lC{$.009$} & \lC{$.17$}  \\ \hline
    $4q-$  &$0$ &$0$ &\lC{$.29$} &$0$ &\lC{$.09$} &\lC{$.05$} &\lC{$.12$}   \\ \hline
    $8p-$  & $0$ & $0$ & \lC{$.07$} & \lC{$.08$} & $0$ & \lC{$.004$} & \mC{$\mathbf{.59}$} \\ \hline
        $1q+$ & $0$ & $0$ & $0$ & \lC{$.004$} & $0$ & $0$ & $0$  \\ \hline
    $Xp-$ & $0$ & $\lC{.001}$ & $\lC{.003}$ & $\lC{.004}$ & \lC{$.01$} & \lC{$.01$} & $0$ \\ \hline
      \end{tabular}

\end{center}
\caption[Estimated confidence for progression model of ovarian
cancer]{{\bf Estimated confidence for ovarian progression.} Frequency
  of edge occurrences in the non-parametric bootstrap test, for the
  trees shown in Figure \ref{fig:treesOvarian}. Colors represent
  confidence: \NewLast{light gray is $< 0.4$, mid gray is in the range
    $[0.4, 0.8]$ and  dark gray is  $> 0.8$.  Bold entries are the edges recovered by the algorithms.}}
\label{fig:freqOvarian}
\end{figure}

\begin{figure}[!htbp]
\begin{center}
\includegraphics[width=0.75\textwidth]{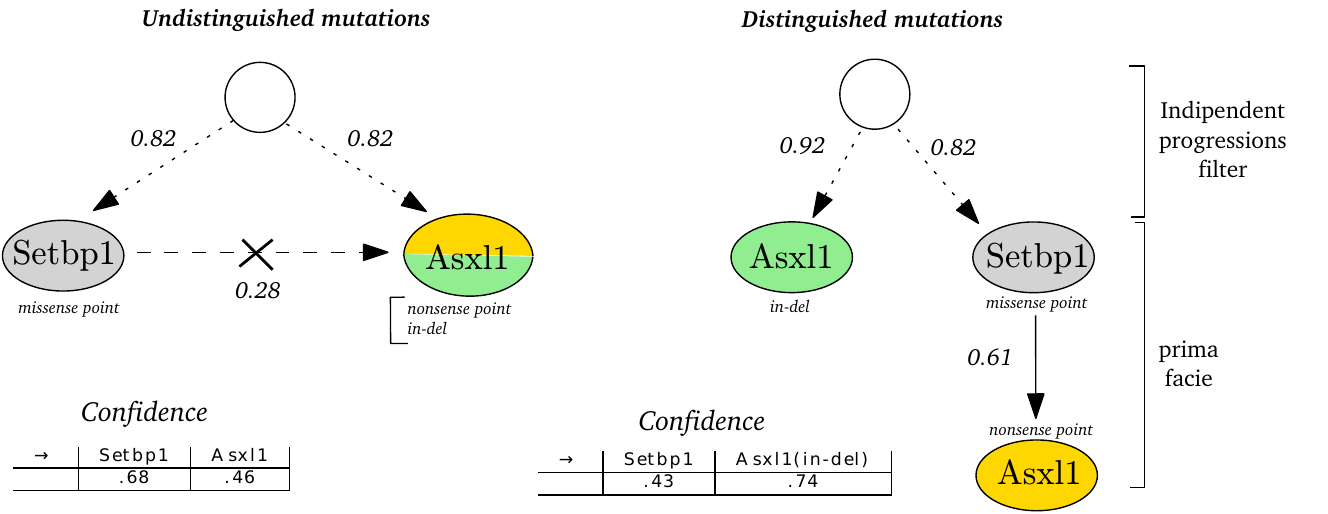}
\end{center}
\caption[Validating the \setbp{}-\asxl{} relation in atypical Chronic Myeloid Leukemia]{{\bf Validating the \setbp{}-\asxl{} relation in atypical Chronic Myeloid Leukemia.} Progression models where \asxl{}  indel and non-sense point are merged (left) and separate (right) suggest that a missense point  mutation  hitting \setbp{} can cause a non-sense point mutation in \asxl{}, that  the observed \asxl{} mutations might be independent and that indel \asxl{} is an early event with high confidence.
}
\label{fig:treesACML}
\end{figure}


\section{Supplementary Materials} 

\subsection{Proofs}\label{sec:proofs}
Here the proofs of all the propositions and theorems follow.

\paragraph{Proof of Proposition \ref{prop:PRdepend} (Dependency).} 
\begin{proof} For $\Rightarrow$ write $\Pconj{\~a}{b}=\Probab{b}- \Pconj{a}{b}$, then write the \PR{} as
\[
\frac{\Pconj{a}{b}}{\Probab{a}} >\frac{\Probab{b} - \Pconj{a}{b}}{1 - \Probab{a}}
\]
and, since $0<\Probab{a}<1$, the proposition follows by simple algebraic arrangements of 
$\Pconj{a}{b} \cdot [1 - \Probab{a}] > \Probab{a}\Probab{b} - \Pconj{a}{b} \cdot \Probab{a}
$. The derivations are analogous but in reverse order for the implication $\Leftarrow$.
\end{proof}

\paragraph{Proof of Proposition \ref{prop:mutualPR} (Mutual probability raising).} 
\begin{proof} The proof follows by Property \ref{prop:PRdepend} and the subsequent implication: 
\[
\Pcond{b}{a} > \Pcond{b}{\~ a} \Leftrightarrow \Pconj{a}{b}>\Probab{a}\Probab{b} \Leftrightarrow
\Pcond{a}{b} > \Pcond{a}{\~ b}\, . 
\]
\end{proof}

\paragraph{Proof of Proposition \ref{prop:PRtemporalPrior} (Probability raising and temporal priority).} 
\begin{proof}
We first prove the forward direction $\Rightarrow$.  Let $x = \Pconj{\~a}{b}$, $y = \Pconj{a}{b}$ and $z =\Pconj{a}{\~b}$.  
We have two assumptions we will use later on: 
\begin{enumerate}
\item $\Probab{a} >\Probab{b}$ which implies $\Pconj{\~a}{b} < \Pconj{a}{\~b}$, i.e., $x < z$. 
\item  $\Pcond{a}{b} > \Pcond{a}{\~ b}$ which, when $0<x+y<1$, implies by simple algebraic rearrangements the inequality 
\begin{equation}
y[1 -  x - y - z] > xz \label{eqn: pr}\, .
\end{equation}
\end{enumerate}
We proceed by rewriting $\Pcond{b}{a}/{\Pcond{b}{\~ a}} >
{\Pcond{a}{b}}/{\Pcond{a}{\~ b}}$ as
\[
\frac{\Pconj{a}{b}\Probab{\~a}}{\Pconj{\~a}{b}\Probab{a}} > 
\frac{\Pconj{a}{b}\Probab{\~b}}{\Pconj{a}{\~b}\Probab{b}}
\]
which means that
\begin{equation}\label{eqn: simplif} 
\frac{\Pcond{b}{a}}{\Pcond{b}{\~ a}} > 
\frac{\Pcond{a}{b}}{\Pcond{a}{\~ b}} 
 \iff 
\frac{\Probab{\~a}}{\Pconj{\~a}{b}\Probab{a}} >
\frac{\Probab{\~b}}{\Pconj{a}{\~b}\Probab{b}} 
\end{equation}

We can rewrite the right-hand side of  (\ref{eqn: simplif}) by using $x$, $y$, $z$
where  $\Probab{a}=\Probab{a,b}+\Probab{a,\~b}=y+z$
and $\Probab{b}=\Probab{a,b}+\Probab{\~a,b}=x+y$, 
 and then do suitable algebraic manipulations. We have
\begin{equation}
\frac{1 - y - z}{x(y +z)}  >  \frac{1 - x - y}{z(x +y)} \iff
yz - y^2z - xz^2 - yz^2 >  xy - x^2y -
x^2z - xy^2 \label{eqn: algeb}
\end{equation}
when $x(y +z)\neq 0$ and $z(x+y)\neq 0$. To check that the right side of  (\ref{eqn: algeb}) holds we show that 
\[
(xy - x^2y - x^2z - xy^2) - (yz - y^2z - xz^2 - yz^2) < 0\, .
\]
First, we rearrange it as $(x-z)[y - y^2 - xz - y(x + z)]<0$  to show that
\begin{align}
(x-z)[y(1- y -x -z) - zx]<  0\label{eqn: end}
\end{align}
is always negative. By observing that, by assumption $1$ we have $z > x$ and thus $(x - z) < 0$, and, 
by equation (\ref{eqn: pr}) we have $y(1- y -x -z) - zx > 0$, we derive
\[
\frac{\Pcond{b}{a}}{\Pcond{b}{\~ a}} >
\frac{\Pcond{a}{b}}{\Pcond{a}{\~ b}} 
\]
which  concludes the $\Rightarrow$ direction. 

The other direction $\Leftarrow$ follows immediately by contraposition: assume 
that $\Pcond{a}{b} > \Pcond{a}{\~ b}$, ${\Pcond{b}{a}}/{\Pcond{b}{\~ a}} > {\Pcond{a}{b}}/{\Pcond{a}{\~ b}}$ 
and $\Probab{b} \leq 
\Probab{a}$.  We distinguish two cases: 
\begin{enumerate}
\item  $\Probab{b} =  
\Probab{a}$, then ${\Pcond{b}{a}}/{\Pcond{b}{\~ a}} =
{\Pcond{a}{b}}/{\Pcond{a}{\~ b}}$.
\item $\Probab{b} <  
\Probab{a}$, then by symmetry  $\Pcond{b}{a} > \Pcond{b}{\~ a}$, and by the $\Rightarrow$ direction of 
the proposition it follows that ${\Pcond{b}{a}}/{\Pcond{b}{\~ a}} <
{\Pcond{a}{b}}/{\Pcond{a}{\~ b}}$.
\end{enumerate}
In both cases we have a contradiction. This completes the proof. 
\end{proof}

\paragraph{Proof of Proposition \ref{prop:monotonicNorm} (Monotonic normalization).} 
\begin{proof}
We prove the forward direction $\Rightarrow$, the converse follows by a similar argument.  Let us assume 
\begin{align}\label{assumption}
\frac{\Pcond{b}{a}}{\Pcond{b}{\~ a}} >
\frac{\Pcond{a}{b}}{\Pcond{a}{\~ b}} 
\end{align}
then $\Pcond{b}{a}\Pcond{a}{\~ b} > \Pcond{a}{b}\Pcond{b}{\~ a}$.
Now, to show the righthand side of the implication, we will show that 
\[
 \Big[\Pcond{b}{a} - \Pcond{b}{\~ a}\Big]\Big[\Pcond{a}{b} + \Pcond{a}{\~ b}\Big] >\Big[\Pcond{b}{a} + 
 \Pcond{b}{\~ a}\Big]\Big[\Pcond{a}{b} - \Pcond{a}{\~ b}\Big]
\]
which reduces to show 
\[
\Pcond{b}{a}\Pcond{a}{\~ b}  - \Pcond{b}{\~ a}\Pcond{a}{b} > \Pcond{b}{\~
  a}\Pcond{a}{b} - \Pcond{b}{a}\Pcond{a}{\~ b} \, .
\]
By (\ref{assumption}), two equivalent inequalities hold
\begin{align*}
 \Pcond{b}{a}\Pcond{a}{\~ b}  - \Pcond{b}{\~ a}\Pcond{a}{b} &> 0\\
\Pcond{b}{\~  a}\Pcond{a}{b} - \Pcond{b}{a}\Pcond{a}{\~ b}  &<0
 \end{align*}
and hence the implication holds. 
\end{proof}

\paragraph{Proof of Proposition \ref{prop:betaDependency} (Coherence in dependency and temporal priority).} 
\begin{proof}
We make two assumptions: 
\begin{enumerate}
\item  $\Pcond{b}{a} > \Pcond{b}{\~ a}$ which implies $\alpha_{a \to b}>0$.
\item $\Probab{a,b} >\Probab{a}\Probab{b}$ which implies $\beta_{a \to b}>0$.
\end{enumerate}
The proof for dependency follows by Proposition \ref{prop:PRdepend} and its implication:
\[
\Pcond{b}{a} > \Pcond{b}{\~ a} \Leftrightarrow \Pconj{a}{b}>\Probab{a}\Probab{b} 
\Leftrightarrow \alpha_{a \to b}>0 \Leftrightarrow \beta_{a \to b}>0. 
\]
Moreover, being $\beta$ symmetric by definition, the proof for temporal priority follows directly by Proposition \ref{prop:mutualPR}
\end{proof}

\New{
The properties outlined so far are sketched informally in the following diagram.

\[
\includegraphics[width=0.9\textwidth]{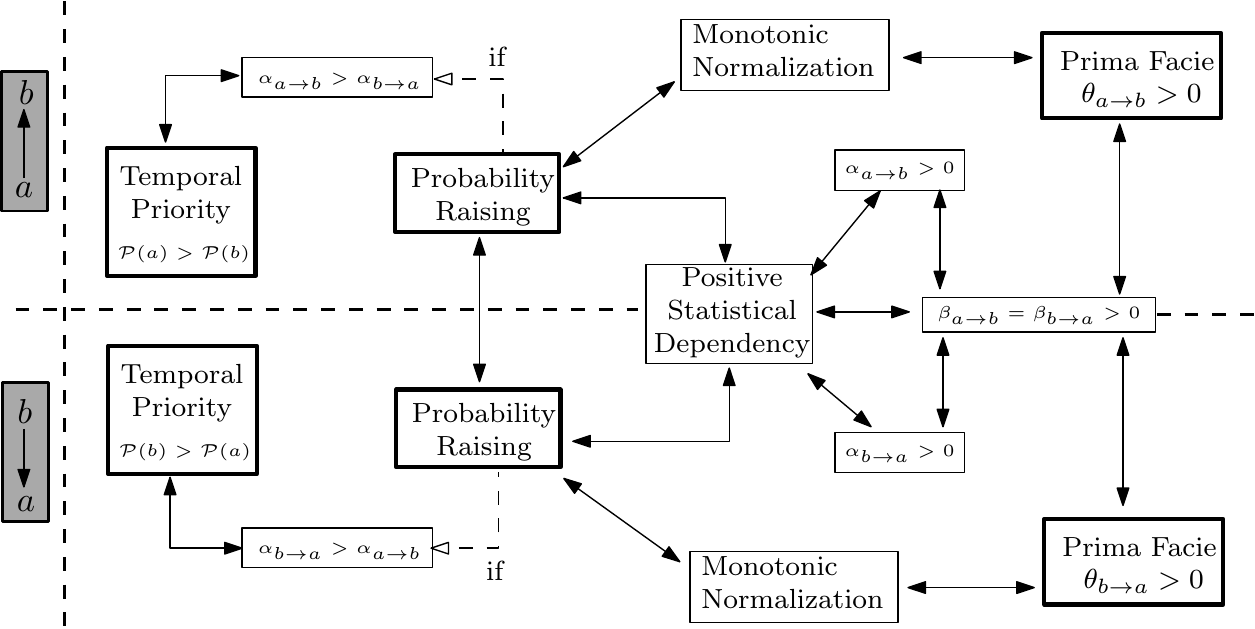}
\]
 }
\paragraph{Proof of Theorem \ref{th:tree-independent} (Independent progressions).}
\begin{proof} For any $a_i \in G^\ast$ it holds that
\[
\Probab{a_i} > \Probab{b} \qquad
m_{a_i\to b} > 0
\]
being prima facie, also $a^\ast \to b$ is the edge selected  by 
\New{CAPRESE} being the $\max\{\cdot\}$ over $G^\ast$.  Thus, 
\New{CAPRESE} selects $\diamond \to b$  instead of $a^{\ast} \to b$  if, for any $a_i$, it holds
\[
\dfrac{1}{1+\Probab{b}} > \dfrac{\Probab{a_i}}{\Probab{a_i}+\Probab{b}}\dfrac{\Pconj{a_i}{b}}{\Probab{a_i}\Probab{b}}\, .
\]
With some algebraic manipulations we   rewrite this as
\[
\Probab{a_i}\Probab{b} + \Probab{b}^2 > \Pconj{a_i}{b}(1+\Probab{b})
\,,\]
which gives the inequality in the theorem statement. 
%
%
\end{proof}

\paragraph{Proof of Theorem \ref{th:tree-nocycles} (Algorithm correctness).}
\begin{proof} It is clear that 
\New{CAPRESE} does not create disconnected components since, to each 
node in $G$, a unique parent is attached (either from $G$ or $\diamond$). For the same reason, no transitive 
connections can appear.

The absence of cycles results from Properties \ref{prop:PRtemporalPrior}, \ref{prop:monotonicNorm} and \ref{prop:betaDependency}. Indeed, 
suppose for contradiction that there is a cycle $(a_1,a_2), (a_2,a_3), 
\ldots, (a_n, a_{1})$ in $E$, then by the three propositions  we have 
\[
\Probab{a_1} >\Probab{a_2}>  \ldots > \Probab{a_n} >  \Probab{a_1} 
\]
which is a contradiction. 
\end{proof}

\paragraph{Proof of Theorem \ref{th:tree-correct} (Asymptotic convergence).}

\begin{proof} 
\New{For any $u \in G$, when $s \to \infty$  the {\em observed probability} $\Probab{u}$  (evaluated from $D$) is equivalent to the product of the probabilities (in $T$) obtained by traversing the forest from the root $\diamond$ to $u$ (Definition \ref{def:treedistrib}). Thus, $\Probab{u} \in (0,1)$ since the traversal probabilities are  in  $(0,1)$ too,  hence all events are distinguishable and  Algorithm 1 reconstructs a tree with the same events set $G$ of $T$.
 
 We now observe that the distribution induced by $T$  (Definition \ref{def:treedistrib}) respects a {\em single-cause prima facie topology} where to each event is assigned {\em at most} a single cause. In other words, Definition \ref{def:praising} holds for any edge $(u,v) \in E$:
 \begin{itemize}
\item by the event-persistence property usually assumed in cancer (fixating mutations are present in  the progeny of a clone)  the occurring times satisfy $t_u < t_v$ which,  in a frequentist sense, implies $\Probab{u} > \Probab{v}$;
\item it holds by construction (Definition \ref{def:treedistrib}) that $\Pconj{v}{u}=\Probab{v}$, thus
$\Pcond{v}{u}=\Probab{v}/{\Probab{u}}$ which is strictly positive since $\Probab{v}$ and $\Probab{u}$ are, and that $\Pconj{v}{\~u}=0$, thus $\Pcond{v}{\~u} = 0$.
\end{itemize}
To correctly reconstruct $T$ we rely on the fact that our score $m_{u\to v}$ is consistent with the prima facie probabilistic causation because of:
 \begin{itemize}
\item  Proposition \ref{prop:PRtemporalPrior}, which states that \PR{} (embedded as $\alpha_{u\to v}$ in $m$) subsumes a good temporal priority model of occurring times, as stated above;
\item Proposition  \ref{prop:monotonicNorm} and \ref{prop:betaDependency} which ensure the monotonicity and sign coherency among $\alpha_{u\to v}$ and $\beta_{u\to v}$ in $m$.
\end{itemize}
Thus, $m$ is consistent with a single-cause prima facie topology. We now show that  Algorithm 1 reconstructs correctly a generic edge in $E$, and hence also $T$.

Consider an event $v \in G$ and  edge $(u,v) \in E$. The set of  its ``candidate'' parent  events  is  $G \setminus \{ v \}$, we  partition it in three disjoint sets $\cal G$, $\cal S$ and $\cal N$:
\begin{itemize}
\item {\em $\cal G$, genuine:}   all the  backward-reachable events,  in $G \setminus \{ v \}$, from $v$;

\item {\em $\cal S$, spurious (or ambiguous):}   all the events (but $v$ ) in the sub-forest which includes the path from $\diamond$ to $v$, which are not in   $\cal G$;

\item {\em $\cal N$, non prima facie:} all other events, i.e., $G \setminus ( \{ v, \} \cup \cal S \cup \cal G)$;
\end{itemize}
Notice that $G = \{ v \} \cup {\cal G} \cup {\cal S} \cup {\cal N}$, that  $u \in \cal G$  and that all the effects of $v$ are non prima facie to $v$ because of the temporal priority, as shown below.}
\[
\includegraphics[width=0.6\textwidth]{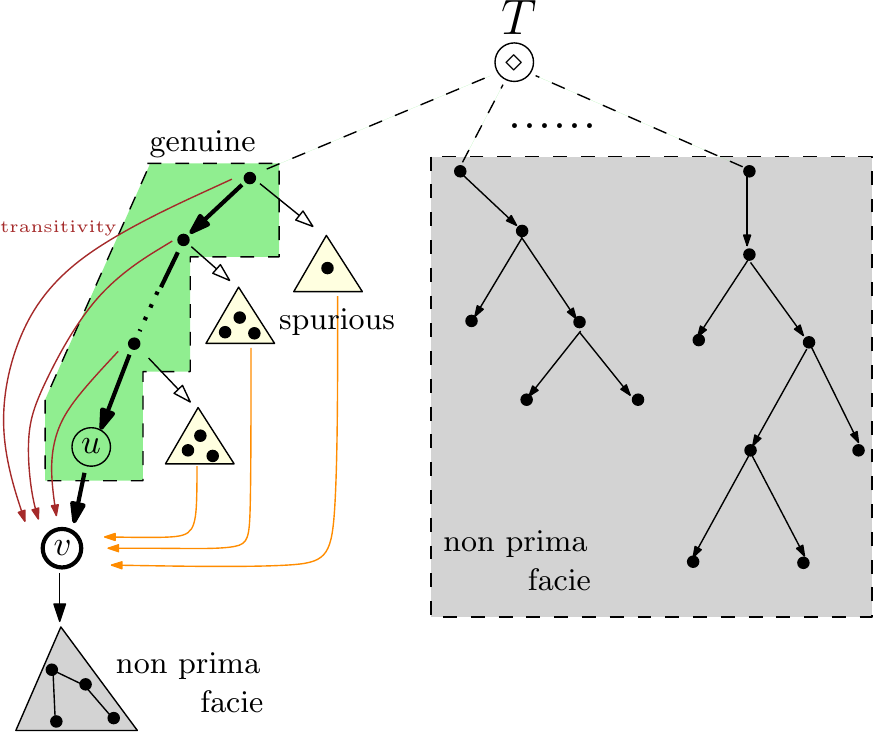}
\]
 
 \New{
This way of partitioning events according to the structure of $T$ subsumes a equivalent partitioning based on the score $\alpha \in [0,1]$, which we use to prove correctness of our algorithm: for any $x \in G$ it holds that $\alpha_{x\to v} = 1$ if $x \in \cal G$,   $0 < \alpha_{x\to v} < 1$ if $x \in \cal S$ and  $\alpha_{x\to v} < 0$ if $x \in \cal N$.

We now show that 
\New{CAPRESE} correctly selects  $u \in \cal G$:
\begin{itemize}
\item  a non prima facie event $x \in \cal N$ \NewLast{either satisfies} (Proposition \ref{prop:PRdepend})
\[
\Pconj{x}{v} \leq \Probab{x} \Probab{v}
\]
which means  that $\alpha_{x \to v} < 0$,  $\beta_{x \to v}  < 0$ (Proposition \ref{prop:betaDependency}) and thus $m_{x \to v} <0$\NewLast{, or it is a descendant of $v$, which means that $\Probab{x}<\Probab{v}$}. By construction, 
\New{CAPRESE} considers as candidate parents of $v$ only  \NewLast{not descendant} events with positive score  (see step 3);

\item a spurious event $x \in \cal S$ is prima facie to $v$ but $\alpha_{x\to v} < 1$ since:
\begin{itemize}
\item $\Probab{x}\Probab{v} < \Pconj{v}{x} < \Probab{v}$, otherwise $x$ would be backward reachable from $x$ and thus in $\cal G$;
\item $0 < \Pconj{v}{\~x} = \Probab{v} - \Pconj{v}{x}$ which means that $\Pconj{v}{\~x}   < \Probab{\~x}\Probab{v}$;
\item by all of the above $\Pcond{v}{\~x} > 0$ which implies that $\alpha_{x \to v} < 1$. 
\end{itemize}
Recall now  that $\lambda \to 0$, which means that $m_{x \to v}  \approx \alpha_{x \to v} <1$. 
\New{CAPRESE} will thus not select any of these events as cause of $v$ if there exist an event with $m_{x \to v} =1$, which is actually the case with genuine causes;
\item genuine causes are the real cause of $v$, $u$, plus {\em all} the transitive  backward-reachable events. Any $x$ of these has maximum score $\alpha_{x\to v} = 1$ since:
\begin{itemize}
\item $\Probab{x}\Probab{v} < \Pconj{v}{x} = \Probab{v}$ and $0 = \Pconj{v}{\~x}$;
\item by the above $\Pcond{v}{\~x} = 0$ which implies  $\alpha_{x \to v} = 1$. 
\end{itemize}
Thus,  
\New{CAPRESE} will pick an event from $\cal G$, and not from $\cal S$. We need to show that $u$ is the event with maximum score.

Enumerate the events in  $\cal G$ as $g_1$ (which is $u$), $\ldots$, $g_k$ in a way that 
\[
\Probab{g_1} <  \ldots < \Probab{g_k}
\] 
and recall that this is a {\em total ordering}  induced by the  temporal priority, and that this is consistent with coefficient $\beta$, which means that 
\[
\beta_{g_1 \to v} >  \ldots > \beta_{g_k \to v}\, .
\] 
Thus, in the limit $\lambda \to 0$ 
\[
\max \{ m_{g_i \to v}  \mid g_i  \in {\cal G}\} \stackrel{\lambda\to 0}{\approx} 
1 + \max \{ \beta_{g_i \to v}  \mid g_i  \in {\cal G}\}  \approx 1+\beta_{g_1 \to v} \approx m_{u \to v} 
\]
is the  event closer in time to $v$, with respect to $\beta$. This
event, namely $u$, is  chosen by the algorithm as the real cause of $v$.
\end{itemize}

Finally, we show that the last step of the algorithm (the independent progression filter, step $4$), does not invalidate the edge $(u,v)$. In fact, the algorithm would  replace such an edge with $(\diamond,v)$ if, for all nodes $x$ backward-reachable from $v$ (i.e., those in ${\cal G} \cup {\cal S}$) it was
\[
\dfrac{1}{1 + \Probab{v}} > \dfrac{\Probab{x}} {\Probab{x} + \Probab{v} } 
 \dfrac{\Pconj{x}{v}} {\Probab{x} \Probab{v} }  \, .
\]
It suffices thus to show that the above inequality is violated just by one of the backward-reachable nodes. We pick just  $u\in {\cal G}$ and note that
\[
 \dfrac{\Probab{u}} {\Probab{u} + \Probab{v} } 
 \dfrac{\Pconj{u}{v}} {\Probab{u} \Probab{v} } 
=
\dfrac{\Pcond{u}{v}} {\Probab{u} + \Probab{v}} \, .
\]
Also, we have that $\Probab{u} < 1$, $\Probab{v} < 1$ and, by
construction, $\Pcond{u}{v}=1$ because all the instances of $v$ are
co-occurring with those of $u$ (but not the converse). Thus, inequality 
\[
\dfrac{1}{1 + \Probab{v}} < \dfrac{1} {\Probab{u} + \Probab{v} }\, ,
\]
is always true and ensures that edge $(u,v)$ is maintained, which concludes the proof. 
}
\end{proof}

\paragraph{Proof of Corollary \ref{cor:noiseunif} (Uniform noise).}

\begin{proof}
\New{As shown in \cite{szabo2002}, the uniform  rates  $\epsilon_+$ and  $\epsilon_-$ affect the observed probabilities as follows
\begin{align}
\Probab{i}^\ast &=\Probab{i} (1-\epsilon_-) + (1-\Probab{i} )\epsilon_+  \\
\Pconj{i}{j}^\ast &=\Pconj{i}{j} (1-\epsilon_-)^2 + [\pi_{ij} - \Pconj{i}{j}](1-\epsilon_-)\epsilon_+ +(1-\pi_{ij})\epsilon_+^2 \, ,
\end{align}
where $\pi_{ij}=\Probab{i} +\Probab{j} -\Pconj{i}{j}$. It is important to note (Lemma 1, \cite{szabo2002}) that
\[
\Probab{i} > \Probab{j} \implies \Probab{i}^\ast > \Probab{j}^\ast \, ,
\]
namely uniform noise  is still implying temporal priority. Because of
this, \NewLast{and since the raw estimate $\alpha$ is monotonic
  relative to temporal priority}, all the derivations for Theorem \ref{th:tree-correct} are still valid in this context, and the algorithm  selects the correct genuine cause for each effect. 

To guarantee that no valid connection is broken by the independent progressions filter, we again rely on Szabo's result (Reconstruction Theorem 1, \cite{szabo2002}). In particular, for any correctly selected edge $(u,v)$ in our algorithm, since we implement Desper's  filter (or, analogously, Szabo's) for independent progressions we do not mistake by deleting $(u,v)$ unless also their algorithms do. Since this is not the case when $\epsilon_+ < \sqrt{p_{\min}}(1-\epsilon_+ - \epsilon_-)$ the proof is concluded.
}
\end{proof}

\subsection{Synthetic data generation}\label{sec:syndatagen}

\New{A set of random trees is generated to prepare synthetic tests. 
Let $n$ be the number of considered events and let $p_{\min}=0.05=1-p_{\max}$,
a {\em single tree} with maximum depth $\log(n)$ is generated as follows:
\begin{algorithmic}[1]
\STATE pick an event $r\in G$ as the tree root;
\STATE  assign to each event but $r$ an integer value in $[2, {\log(n)}]$ representing its depth in the tree,
 ensure that for each level there is at least one event ($0$ is reserved for $\diamond$, $1$ for $r$);
 \FORALL{events $e \neq r$} \STATE{let $l$ be the level assigned to $e$;}
 \STATE{assign a father to $e$ selecting an event among those at which level $l-1$ was assigned;} 
 \STATE{add the selected pair to the set of edges $E$}; 
 \ENDFOR
 \FORALL{edges $(i,j)\in E$} 
  \STATE{assign $\alpha((i,j))$ a random value in $[p_{\min},p_{\max}]$;} 
 \ENDFOR
\RETURN the generated tree;
\end{algorithmic}
When a forest is to be generated, we repeat the above algorithm to create its \NewLast{constituent} trees. These trees (or forests), in turn, are used to sample the input matrix for the reconstruction algorithms, with the  parameters  described in the main text. 
}

\subsection{Further results} \label{sec:figures}

We show here the results of the experiments discussed but not presented in the main text.

\begin{figure}[t]
\begin{center}
\includegraphics[width=0.99\textwidth]{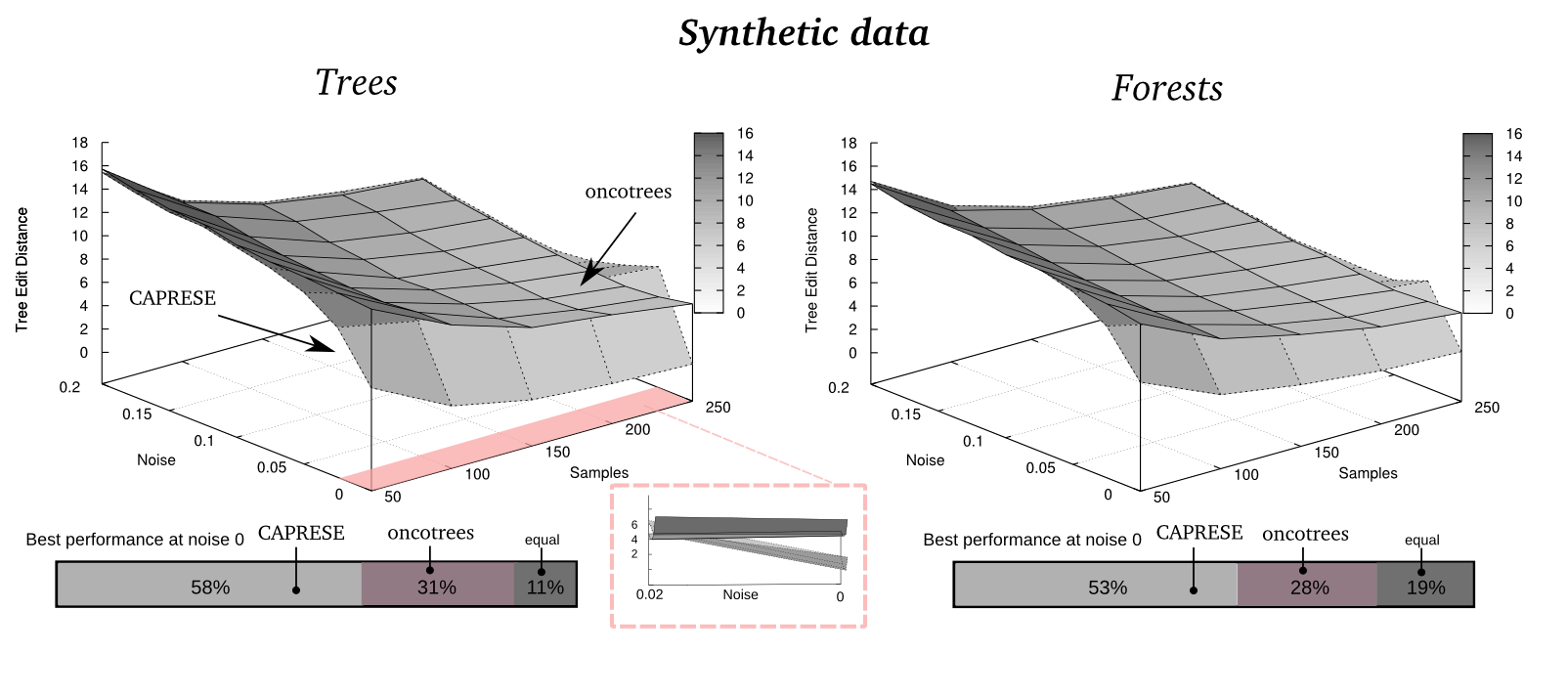}
\end{center}
\caption[Reconstruction with noisy synthetic data and $\lambda \to 0$.]{{\bf Reconstruction with noisy synthetic data and $\lambda \to 0$.} The settings of the experiments are the 
same as those used in  Figure \ref{fig:synteticnoisy}, but in this case the estimator is shrank by $\lambda \to 0$, i.e., $\lambda=0.01$. In the magnified image one can sees that the performance of 
\New{CAPRESE} converges to Desper's one already for $\nu \approx 0.01$, hence largely faster than in the case of $\lambda \approx 1/2$ (Fig. \ref{fig:synteticnoisy}). }
\label{fig:synteticnoisyshrinkalpha}
\end{figure}

\paragraph{Reconstruction of noisy synthetic data with $\lambda \to 0$.}

Although we know that $\lambda \to 0$ is not the optimal value of the shrinkage-like coefficient for noisy data, we show in Figure \ref{fig:synteticnoisyshrinkalpha} the analogue of Figure \ref{fig:synteticnoisy} when the estimator is shrank by $\lambda \to 0$, i.e., $\lambda=0.01$. When compared  to Figure \ref{fig:synteticnoisy} it is clear that a best performance of 
\New{CAPRESE}  is obtained with $\lambda \approx 1/2$, as suggested by Figure \ref{fig:lambda}.

\begin{figure}[!t]
\begin{center}
\includegraphics[width=0.60\textwidth]{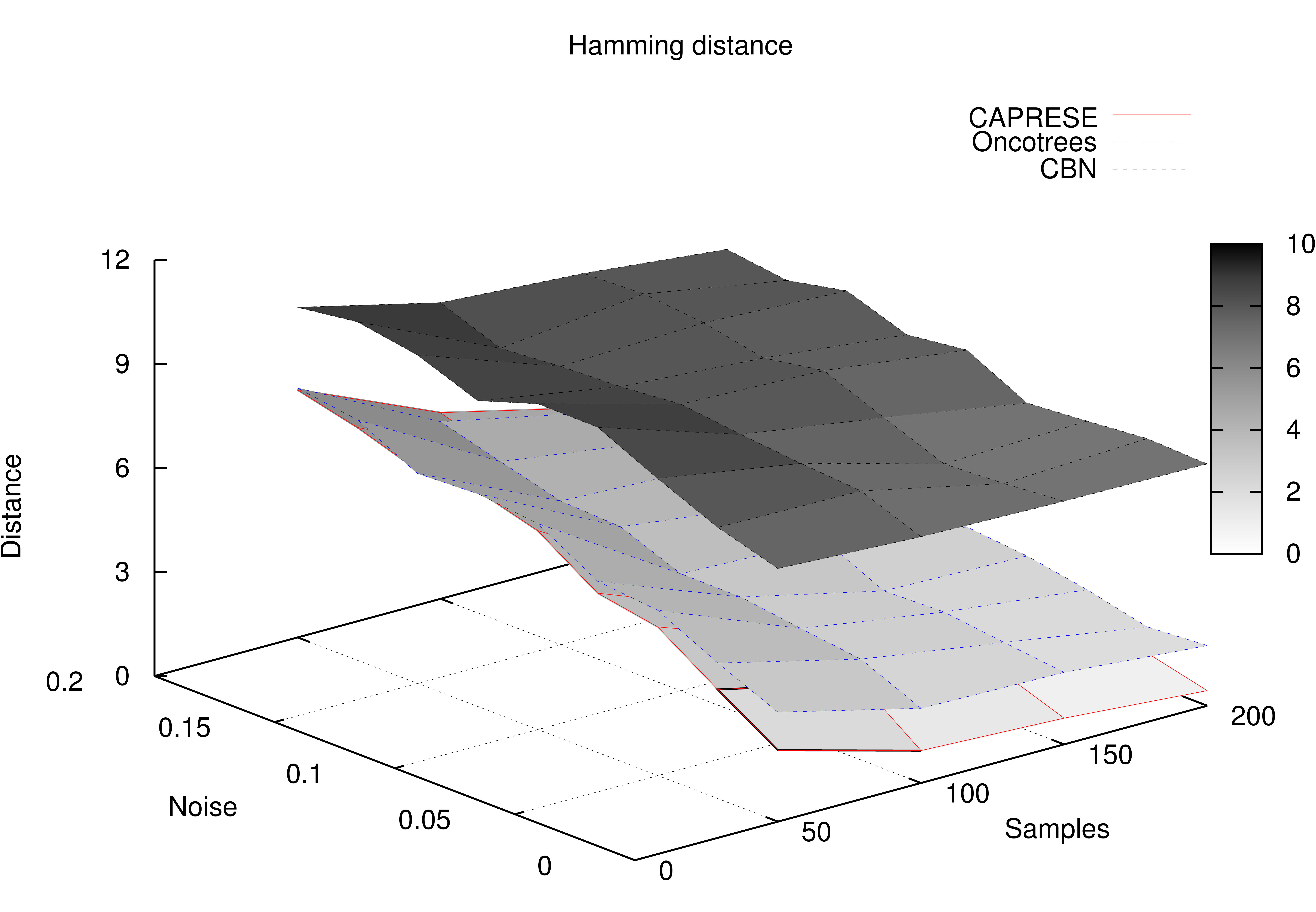}
\end{center}
\caption[Reconstruction with {CAPRESE} compared to oncotrees and h-CBN  with noisy synthetic data.]{{\bf Reconstruction with {CAPRESE} compared to oncotrees and h-CBN  with noisy synthetic data.} Performance of 
{CAPRESE} compared to oncotrees and h-CBNs  as a
function of the number of samples and noise $\nu$. The $\lambda$
parameter used for 
\New{CAPRESE} is 
$1/2$, and the reconstructed topologies contain 10 nodes each.}
\label{fig:synteticnoisycbn}
\end{figure}

\New{
\paragraph{Comparison with hidden Conjunctive Bayesian Networks, h-CBNs.}
We here compare the performance of CAPRESE to hidden Conjunctive Bayesian
Networks (h-CBN)\cite{gerstung2009quantifying}, as well as to oncotrees. The settings of the experiment are slightly different from those of the previous analyses: we used
100 distinct random trees of 10 events each. We ranged the number of samples available
for reconstruction from 50 to 200, with a step size of 50. The settings
used for running h-CBNs are relatively standard settings: we allowed
for 50 annealing iterations with initial temperature equal to 1. 
Since h-CBNs reconstruct DAGs, it is not possible to quantify its
performance using Tree Edit Distance, as we did in the comparison with oncotrees. Instead, we here adopt Hamming Distance (computed on the connection adjacency matrix), as a closely related and computationally feasible alternative for measuring performance \cite{hamming_distance}. 

The results of the experiment can be found in Figure \ref{fig:synteticnoisycbn}, and
show  that CAPRESE clearly outperforms h-CBNs.  In particular, it is
possible to notice that, for all the analyzed values of noise and
sample sizes, both CAPRESE and
oncotrees display a (average) Hamming Distance between the
reconstructed model and the original tree topology that is
significantly lower than h-CBNs, with the largest differences observed
in the noise-free case. This result would point at a much faster
convergence of CAPRESE with respect to the number of samples, also in presence of moderate levels of noise.  
}

\New{

A few remarks are warranted about this experiment. First, in contrast to the comparison with oncotrees, we ran each experiment exactly once rather
than averaging the results over 10 repetitions, and on relatively smaller trees. These limitations are a
consequence of
the extremely high time
complexity of the simulated annealing step of h-CBNs. However, the
comparison between CAPRESE and h-CBNs shows a so large difference in the performance that we do not
expect this to be have significant impact. Second, the results
obtained by h-CBNs are perhaps worse than expected based on results in
the absence of noise presented in \cite{Hainke:2012hwa}, which were however based on a unique tree topology. Yet, this outcome may have been potentially
influenced by either the estimation procedure of the noise
parameter in h-CBN, the adopted annealing procedure or by the
used number of iterations. 
In future work we plan to extend our algorithm to extract more general topologies and to compare both methods in a greater detail. 
}

\begin{figure}[t]
\begin{center}
\includegraphics[width=0.70\textwidth]{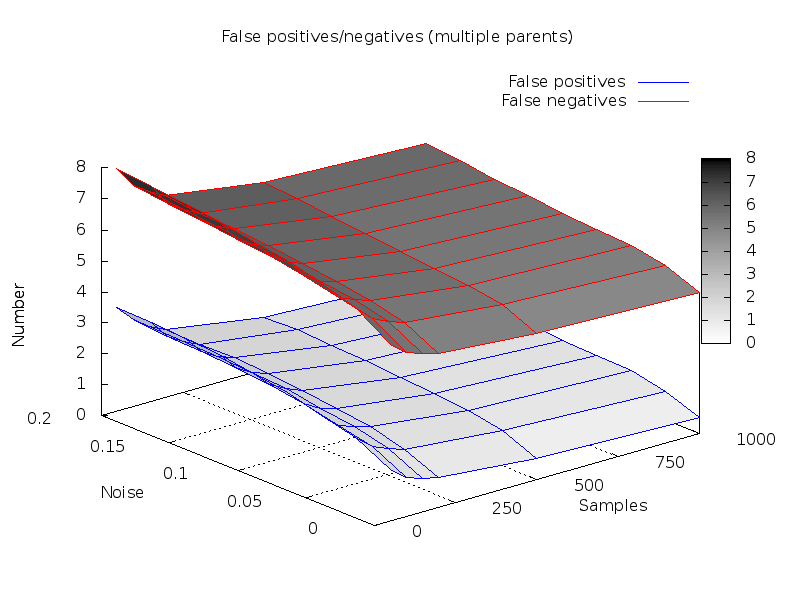}
\end{center}
\caption[Performance of CAPRESE to reconstruct models with conjunctive parents and noisy data.]{
{\bf Performance of CAPRESE to reconstruct models with conjunctive
  parents and noisy data.} \NewLast{Performance of CAPRESE measured in
  terms of the number of {\em false positives/negatives} in the reconstructed model, when
  data are generated from  directed acyclic graphs with $10$ nodes and where each event is caused by at most $3$ {\em conjunctive events} (randomly assigned). The  $\lambda$ parameter is set to $1/2$.}}
\label{fig:DAGs-CAPRESE}
\end{figure}

\NewLast{
\paragraph{Inference of models with multiple conjunctive parents.}
CAPRESE is specifically tailored to reconstruct models with independent progressions and a unique cause for each event (i.e., trees or forests), while other approaches such as CBNs can reconstruct models where multiple conjunctive parents co-occur to cause an effect (i.e., $a \wedge b$ cause $c$). It is thus reasonable to use such conjunctive approaches to infer more complex model, in spite of CAPRESE.

However, it is interesting to asses CAPRESE's performance when
(synthetic) data are sampled from a model with multiple parents and
noise. By sampling input data from random {\em directed acyclic graphs}
with $10$ nodes and where each event is caused by at most $3$
{conjunctive events} (randomly assigned), we assess the number of {\em
  false positives} and {\em false negatives} retrieved in the model reconstructed with
CAPRESE. We show the results in Figure \ref{fig:DAGs-CAPRESE}. Our
results indicate that for increasing sample size, the number of false
positives approaches $0$. Thus, for sufficiently large number of samples,  all the causal claims returned by
CAPRESE are true. In addition, the
number of false negatives is always higher and proportional to the
connectivity of the target model. This is to be expected since CAPRESE
assigns at most one parent (the cause) to every node. 
}

\end{document}